\documentclass{article}
\usepackage{amssymb}
\usepackage{amsmath}
\usepackage{amsfonts}
\usepackage{graphicx}
\newtheorem{thm}{Theorem}[section]
\newtheorem{theorem}{Theorem}

\newtheorem{corollary}[theorem]{Corollary}

\newtheorem{definition}[theorem]{Definition}
\newtheorem{example}[theorem]{Example}

\newtheorem{lemma}[theorem]{Lemma}
\newtheorem{notation}[theorem]{Notation}

\newtheorem{proposition}[theorem]{Proposition}
\newtheorem{remark}[theorem]{Remark}

\newenvironment{proof}[1][Proof]{\noindent\textbf{#1.} }{\ \rule{0.5em}{0.5em}}
\DeclareMathOperator{\aut}{Aut}
\DeclareMathOperator{\Ker}{Ker}
\begin{document}

\title{Finite abelian subgroups of the mapping class group}
\author{S. Allen Broughton, Rose-Hulman Institute of Technology
\and Aaron Wootton, University of Portland }
\maketitle

\begin{abstract}
The problem of enumeration of conjugacy classes of finite abelian subgroups of the mapping class group $\mathcal{M}_{\sigma}$ of a compact closed surface $X$ of genus $\sigma$ is considered. A complete method of enumeration is achieved for finite elementary abelian subgroups and steps are taken toward enumeration of finite abelian subgroups.

\end{abstract}


\begin{center}



\end{center}

\section{Introduction\label{Sec-Intro}}

Let $S$ be a closed, smooth, orientable surface of genus $\sigma\geq2.$ The
mapping class group $\mathcal{M}_{\sigma}$ of $S$ (or MCG) is the group of
isotopy classes of homeomorphisms of $S.$ In this paper we shall investigate
the conjugacy classes of finite abelian subgroups of $\mathcal{M}_{\sigma}$.
While the general finite subgroup classification is important we focus on the
elementary abelian case as a tractable case where complete classification
methods by standard linear algebra is possible and the abelian case where
positive steps can be made toward such a classification. The main result of
our work is to describe methods which may be employed to completely classify
all elementary abelian actions in a given genus and steps toward a
classification of abelian actions. The classification is complex since it
involves understanding the representation theory of certain subgroups of
symmetric groups, but for a fixed low genus since the symmetric groups are
small, we are able to produce very explicit results. \newline

Our main results are presented in three parts. Theorem \ref{th-split} gives a
decomposition of an abelian subgroup $G$ into an unramified part (hyperbolic)
and a ramified part (elliptic). For elementary abelian subgroups, the
classification of the unramified part is completely described in Section
\ref{Sec-Unramified} and methods allowing one to classify the ramified part,
as well as some explicit examples, are in Section \ref{Sec-Ramified}. In
Section \ref{Sec-Unramified}, steps are taken to classify general unramified
abelian actions. Though complete results are not obtained, enough information
is gathered to classify unramified abelian actions on surfaces up to genus
$65$ (see Example \ref{ex-genus65}) and the general method of how to classify
abelian actions is described through other examples. \newline

We should note that our methods fail for non-abelian groups. There is a
general classification theory which describes the classes of subgroups as a
finite sequence of quotients of finite sets by the action of infinite groups.
In the abelian case the general method simplifies greatly and we get a nice
splitting into ramified and unramified cases and the subcases are tractable.
For the general case however, these simplifications do not occur and
calculations become complicated very quickly and are only really possible
through computer calculation. For low genus, some partial results are known,
and many cases of infinite families, as mentioned further in this
introduction. \newline

The importance of a good understanding of the finite subgroups of $M_{\sigma}$
was made apparent when Maclachlan showed that $M_{\sigma}$ is generated by
torsion elements in \cite{Mac2}. Following this, a number of different people
have found sets of torsion generators for $M_{\sigma}$, see \cite{Bi} and
\cite{BreFarb}. An understanding of how such elements generate $M_{\sigma}$
may be possible through an analysis of the conjugacy classes of finite
subgroups of $M_{\sigma}$. Another motivating reason for our work is that each
of the conjugacy classes of finite subgroups of $\mathcal{M}_{\sigma}$ are in
1-1 correspondence to finite groups of orientation preserving groups of
homeomorphisms of $S$ up to topological equivalence. See \cite{Bro2} for a
detailed discussion and the complete classification in genus 2 and 3 (with one
omission \cite{Bro3}). A brief summary of the equivalence is given in
Subsection \ref{subsec-finite} below. Another motivating reason is that the finite
subgroups describe the singularity structure of moduli space with implications
about the structure of the cohomology of the mapping class group. The moduli
space $\mathfrak{M}_{\sigma}$ of surfaces of genus $\sigma$ is the space of
conformal equivalence classes of surfaces of genus $\sigma.$ The moduli space
may be obtained as a quotient of the Teichm\"{u}ller space $\mathfrak{M}%
_{\sigma}=\mathfrak{T}_{\sigma}/\mathcal{M}_{\sigma},$ where $\mathfrak{T}%
_{\sigma}$ is homeomorphic to an open ball in $\mathbb{C}^{3\sigma-3}.$ The
singular points of $\mathfrak{M}_{\sigma}$ are caused by fixed points of
finite subgroups of $\mathcal{M}_{\sigma}$ which in turn controls the torsion
cohomology of the mapping class group. For more details and applications see
\cite{Bro1}. \newline


A starting point for our work are the papers \cite{Har1} and \cite{Mac1}
discussing cyclic and abelian groups of surface automorphisms. We adapt and
extend the methods of \cite{Har2} especially in the discussion of unramified
actions using the correspondence between groups of conformal automorphisms and
conjugacy classes of the MCG given by the Nielsen Realization Problem, see
\cite{K}. The advantage of considering $G\leq M_{\sigma}$ in this way is that
there has been tremendous progress in classification results of automorphism
groups of compact Riemann surfaces. A review of some results are given in the
paper \cite{Bro2} and the monograph \cite{Breu}. Some additional calculations
are available in preprint form \cite{Broet1} and \cite{JGIA}. In addition,
there is much literature dedicated to infinite families of surfaces whose
automorphism groups share certain properties, see for example \cite{Har1} or
\cite{Woo}. For reference, a comprehensive study in this area is given in
Breuer's monograph. We note that in Breuer's monograph, group actions are considered
equivalent if their representations on the first homology group are equivalent - a coarser classification than the
classification up to conjugacy in the mapping class group.

\textbf{Acknowledgements }The first author wishes to thank I.M. Isaacs for
advice shortening the exposition in Section 4. The second author thanks the
University of Portland for support through the Butine faculty support program.

\section{Reduction of classification to finite group
calculations\label{Sec-Reduction}}

It is standard to study surface groups, Teichm\"{u}ller space, and the moduli
space in terms of Fuchsian groups so we present our problem in that context.

\subsection{Finite group actions and finite subgroups of $\mathcal{M}_{\sigma
}$\label{subsec-finite}}

Let $G$ be a finite group. $\ $The group $G$ is said to act (in an orientation
preserving manner) on surface a $S$ of genus $\sigma\geq2$ if there is an
injection
\[
\epsilon:G\hookrightarrow\mathrm{Homeo}^{+}(S)
\]
into the group of orientation preserving homeomorphisms. Two actions
$\epsilon_{1},\epsilon_{2}$ are said to be \emph{topologically equivalent} it
there is a homeomorphism $h$ of \ $S$ and an automorphism $\omega$ of $G$ such
that%
\[
\epsilon_{2}(\omega(g))=h\circ\epsilon_{1}(g)\circ h^{-1}.
\]
This is equivalent to saying that the images $\epsilon_{1}(G)$ and
$\epsilon_{2}(G)$ are conjugate in $\mathrm{Homeo}^{+}(S)$. The following are
well known for $\sigma\geq2$:

\begin{itemize}
\item The map to the mapping class group $\mathcal{M}_{\sigma}=$
$\mathcal{M(}S\mathcal{)}$ is injective
\begin{equation}
G\overset{\epsilon}{\longrightarrow}\mathrm{Homeo}^{+}(S)\overset{\iota
}{\longrightarrow}\mathcal{M}_{\sigma} \label{eq-GtoM}%
\end{equation}
because $G$ acts faithfully on $H_{1}(S).$

\item Given a finite subgroup $H\subset\mathcal{M}_{\sigma}$ there is a finite
group $G\subset$ $\mathrm{Homeo}^{+}(S)$ such that $H=\iota(G)$ and at least
one conformal structure on $S$ such that $G$ is a group of conformal
automorphisms with respect to this conformal structure (Nielsen Realization
problem \cite{K}).
\end{itemize}

Thus we have the following:

\begin{proposition}
The map induced by (\ref{eq-GtoM}) is a 1-1 correspondence between topological
equivalence classes of $G$-actions and conjugacy classes of the mapping class group.
\end{proposition}

\subsection{Covering Fuchsian group and generating vectors}

Now suppose that $G$ acts conformally on $S.$ The quotient $T=S/G$ is a
surface of some genus $\rho$ and the quotient map $S\rightarrow S/G$ is
branched over $r$ points $Q_{1},\ldots,Q_{r}$ with periods (or branching
orders) $m_{1},\ldots,m_{r}.$ We say that $\mathcal{S}=(\rho:m_{1}%
,\ldots,m_{r})$ is the \emph{signature} or \emph{branching data} of $G$ acting
on $S $ and that $S$ admits a $G$-$(\rho:m_{1},\ldots,m_{r})$ action.
\emph{\ }Since $G$ acts conformally there is a uniformizing Fuchsian group
$\Gamma$ for the $G$-action and the signature of $\Gamma$ is $(\rho
:m_{1},\ldots,m_{r}).$ More precisely we have a map
\[
\Pi\hookrightarrow\Gamma\overset{\eta}{\twoheadrightarrow}G,
\]
where $\Pi$ is a torsion free surface group acting on $\mathbb{H}$ so that
$S\approx\mathbb{H}/\Pi,$ and the $G$-action on $S$ is induced by the
isomorphism $\overline{\eta}:\Gamma/\Pi\longleftrightarrow G$ and the natural
action of $\Gamma/\Pi$ on $\mathbb{H}/\Pi.$

\bigskip The epimorphisms $\Gamma\overset{\eta}{\twoheadrightarrow}G$ can be
neatly summarized in the context of finite groups by generating vectors (see
\cite{Bro1} and \cite{Bro2}). For the conformal action of $G$ on a surface $S$
with signature $\mathcal{S}=(\rho:m_{1},\ldots,m_{r}),$ the group $\Gamma$ has
signature $(\rho:m_{1},\ldots,m_{r})$ and a presentation of the following
type:
\begin{align}
\Gamma &  \simeq\Gamma(\rho:m_{1},\ldots,m_{r})\label{eq-Gammapres}\\
&  =\langle\alpha_{i},\beta_{i},\gamma_{j},1\leq i\leq\rho,1\leq j\leq
r:\prod_{i=1}^{\rho}[\alpha_{i},\beta_{i}]\prod_{j=1}^{r}\gamma_{j}=\gamma
_{1}^{m_{1}}=\cdots=\gamma_{r}^{m_{r}}=1\rangle,\nonumber
\end{align}
and the genus $\sigma$ of $S$ is given by the Riemann-Hurwitz equation
\begin{equation}
\frac{(2\sigma-2)}{|G|}=(2\rho-2+r)-\sum_{j=1}^{r}\frac{1}{m_{j}}.
\label{eq-RH}%
\end{equation}
Let $a_{i},b_{i},c_{j}$ be images of the generators $\alpha_{i},\beta
_{i},\gamma_{j}$ under the epimorphism $\eta:$ $\Gamma\longrightarrow G$,
i.e.,
\begin{equation}
\eta:\alpha_{i}\rightarrow a_{i},\text{ }\eta:\beta_{i}\rightarrow
b_{i},\text{ }\eta:\gamma_{j}\rightarrow c_{j}, \label{eq-GVdef}%
\end{equation}
Then the set $\{a_{i},b_{i},c_{j}:1\leq i\leq\rho,1\leq j\leq r\}$ is a
generating set for $G$ satisfying the following properties:
\begin{equation}
\prod_{i=1}^{\rho}[a_{i},b_{i}]\prod_{j=1}^{r}c_{j}=1, \label{eq-GVprop1}%
\end{equation}
and
\begin{equation}
o(c_{j})=m_{j}. \label{eq-GVprop2}%
\end{equation}
Each such generating set gives us a $\left(  2\rho+r\right)  $-tuple
$(a_{1},\ldots a_{\rho},b_{1},\ldots b_{\rho},c_{1},\ldots c_{\rho})$
satisfying (\ref{eq-GVprop1}) and (\ref{eq-GVprop2}) and is called a
$(\rho:m_{1},\ldots,m_{r})$ - \emph{generating vector}. There is a 1-1
correspondence between the set of $(\rho:m_{1},\ldots,m_{r})$ - generating
vectors of $G$ and $Epi(\Gamma,G),$ epimorphisms $\Gamma\rightarrow G$
preserving the orders of the $\gamma_{j},$ once a generating set
$\mathcal{G}=\{\alpha_{i},\beta_{i},\gamma_{j}:1\leq i\leq\rho,1\leq j\leq
r\}$ of $\Gamma$ has been fixed. Let $X^{\circ}(G,\mathcal{S)}=X^{\circ
}(G,\rho:m_{1},\ldots,m_{r})$ denote the set of $(\rho:m_{1},\ldots,m_{r})$ -
generating vectors of $G$ acting on $S$.

There is a natural action of $\mathrm{Aut}(G)\times\mathrm{Aut}(\Gamma)$ on
$Epi(\Gamma,G)$ by
\[
\eta\longrightarrow\omega\circ\eta\circ\xi^{-1}%
\]
where $(\omega,\xi)\in\mathrm{Aut}(G)\times\mathrm{Aut}(\Gamma).$ The action
transfers to $X^{\circ}(G,\mathcal{S)}$ in a natural way. The following is
well known and the justification provided in some detail in \cite{Bro2}.

\begin{proposition}
Let notation be as above. Then we have:

\begin{enumerate}
\item Each finite subgroup of the mapping class group has a uniquely
determined signature which may be recovered from the homology representation.

\item The conjugacy classes of finite subgroups of the mapping class group are
in 1-1 correspondence to the $\mathrm{Aut}(G)\times\mathrm{Aut}(\Gamma) $
orbits on $X^{\circ}(G,\mathcal{S)}.$

\item The group $\mathrm{Aut}(\Gamma)$ may be constructed geometrically from
the mapping class group $\mathcal{M}_{\mathcal{B}}$ of $T$ preserving the
branch point set and the branching orders.
\end{enumerate}
\end{proposition}

\begin{remark}
\label{rk-FundG}Let's give a bit more detail on the last point of the
proposition. Let $\mathcal{B}=\{Q_{1},\ldots,Q_{r}\}$ denote the branch point
set, and $T^{o}=T-\mathcal{B}$ and $Q_{0}\in T^{o}.$ Every homeomorphism of
$T^{o}$ may be identified with a homeomorphism of the pair $(T,\mathcal{B}%
)\ $and vice versa. The fundamental group, $\pi_{1}(T^{o},Q_{0})$ has a
presentation of the form:%
\begin{equation}
\langle\widetilde{\alpha}_{i},\widetilde{\beta}_{i},\widetilde{\gamma}%
_{j},1\leq i\leq\rho,1\leq j\leq r:\prod_{i=1}^{\rho}[\widetilde{\alpha}%
_{i},\widetilde{\beta}_{i}]\prod_{j=1}^{r}\widetilde{\gamma}_{j}%
=1\rangle\label{eq-FundGpres}%
\end{equation}
where the $\widetilde{\alpha}_{i},\widetilde{\beta}_{i}$ are the canonical
pair around the $i$'th handle and $\widetilde{\gamma}_{j}$ encircles the
$j$'th puncture. There is an epimorphism $\pi_{1}(T^{o},Q_{0})\rightarrow
\Gamma$ given by%
\begin{equation}
\widetilde{\alpha}_{i}\rightarrow\alpha_{i},\widetilde{\beta}_{i}%
\rightarrow\beta_{i},\widetilde{\gamma}_{j}\rightarrow\gamma_{j}.
\label{eq-FundtoGamma}%
\end{equation}
Let $\mathcal{M(}T^{o},\mathcal{S)}$ denote the subgroup of the mapping class
group of $\mathcal{M(}T^{o}\mathcal{)}$ that preserves periods of the branch
points. Then the canonical map $\mathcal{M(}T^{o},\mathcal{S)\rightarrow
}\mathrm{Out}(\Gamma)$ is an isomorphism, see \cite{Z}.
\end{remark}

\subsection{The $\mathrm{Aut}(G)\times\mathrm{Aut}(\Gamma)$ actions for
abelian groups.\label{subsec-AGxAGam}}

It is clear how $\mathrm{Aut}(G)$ acts on $X^{\circ}(G,\mathcal{S)}$ namely
$\omega\in\mathrm{Aut}(G)$ acts via:
\[
(a_{1},\ldots a_{\rho},b_{1},\ldots b_{\rho},c_{1},\ldots c_{\rho}%
)\rightarrow(\omega a_{1},\ldots\omega a_{\rho},\omega b_{1},\ldots\omega
b_{\rho},\omega c_{1},\ldots\omega c_{r}).
\]
In order to compute the action of $\mathrm{Aut}(\Gamma)$ on $X^{\circ
}(G,\mathcal{S)}$ we need a generating set for $\mathrm{Aut}(\Gamma).$ This
can be constructed as follows. In \cite{Bi} Birman constructs a set of
geometrically defined generators of $\mathcal{M(}T^{o},\mathcal{S)}$ using
Dehn twists and spin maps. The action of these geometric generators on
$\pi_{1}(T^{o},Q_{0})$ can be written as substitution formulas in the
generators $\langle\widetilde{\alpha}_{i},\widetilde{\beta}_{i},\widetilde
{\gamma}_{j},1\leq i\leq\rho,1\leq j\leq r\rangle$ using the presentation
(\ref{eq-FundGpres}) and then writing $\langle\alpha_{i},\beta_{i},\gamma
_{j},1\leq i\leq\rho,1\leq j\leq r\rangle$ as substitution formulas in the
generators using the homomorphism (\ref{eq-FundtoGamma}) and the presentation
(\ref{eq-Gammapres}). Any substitution formula can be tested against the
presentation to see if it really is an automorphism. Here are some examples,
there are many more but we will not record them here since we are only going
to use the much simpler transformations for Abelian groups.
\[%
\begin{tabular}
[c]{c}%
Action of automorphisms $\mathfrak{A}_{i},$ $\mathfrak{B}_{i},$ $\mathfrak{C}%
_{i},\mathfrak{R}_{i},$ $\mathfrak{S}_{i},$ $\mathfrak{Z}_{i}$\\
$%
\begin{tabular}
[c]{|c|c|c|c|c|}\hline
& $\alpha_{i}$ & $\beta_{i}$ & $\alpha_{i+1}$ & $\beta_{i+1}$\\\hline
$\mathfrak{A}_{i}$ & $\alpha_{i}$ & $\beta_{i}\alpha_{i}$ & $\alpha_{i+1}$ &
$\beta_{i+1}$\\\hline
$\mathfrak{B}_{i}$ & $\alpha_{i}\beta_{i}$ & $\beta_{i}$ & $\alpha_{i+1}$ &
$\beta_{i+1}$\\\hline
$\mathfrak{R}_{i}$ & $\alpha_{i}\beta_{i}\alpha_{i}^{-1}$ & $\alpha_{i}^{-1}$
& $\alpha_{i+1}$ & $\beta_{i+1}$\\\hline
$\mathfrak{S}_{i}$ & $\delta\alpha_{i+1}\delta^{-1}$ & $\delta\beta
_{i+1}\delta^{-1}$ & $\alpha_{i}$ & $\beta_{i}$\\\hline
$\mathfrak{Z}_{i}$ & $\alpha_{i}\alpha_{i+1}$ & $\alpha_{i+1}^{-1}\beta
_{i}\alpha_{i+1}$ & $\epsilon\alpha_{i+1}$ & $\beta_{i+1}\beta_{i}%
^{-1}\epsilon^{-1}$\\\hline
\end{tabular}
$\\
\\
$\delta=[\alpha_{i},\beta_{i}],$ $\varepsilon=[\alpha_{i+1}^{-1},\beta_{i}]$%
\end{tabular}
\]

\[%
\begin{tabular}
[c]{c}%
Action of $\mathfrak{T}_{j}$\\
$%
\begin{tabular}
[c]{|c|c|c|c|}\hline
& $\gamma_{j}$ & $\gamma_{j+1}$ & notes\\\hline
$\mathfrak{T}_{j}$ & $\gamma_{j}\gamma_{j+1}\gamma_{j}^{-1}$ & $\gamma_{j}$ &
$o(\gamma_{j})=o(\gamma_{j+1})$\\\hline
$\mathfrak{T}_{j}^{-1}$ & $\gamma_{j+1}$ & $\gamma_{j+1}^{-1}\gamma_{j}%
\gamma_{j+1}$ & $o(\gamma_{j})=o(\gamma_{j+1})$\\\hline
\end{tabular}
$%
\end{tabular}
\]

\[%
\begin{tabular}
[c]{c}%
Action of $\mathfrak{U}_{i,j},\mathfrak{V}_{i,j}$\\
$%
\begin{tabular}
[c]{|c|c|c|c|c|}\hline
& $\alpha_{i}$ & $\beta_{i}$ & $\gamma_{j}$ & notes\\\hline
$\mathfrak{U}_{i,j}$ & $\alpha_{i}$ & $\beta_{i}x\gamma_{j}x^{-1}\ $ &
$y\gamma_{j}y^{-1}$ & $x=\alpha_{i}^{-1}\beta_{i}^{-1}w_{2},$ $y=x^{-1}%
\alpha_{i}^{-1}x$\\\hline
$\mathfrak{V}_{i,j}$ & $\alpha_{i}u\gamma_{j}u^{-1}$ & $\beta_{i}$ &
$v\gamma_{j}v^{-1}$ & $u=\beta_{i}\alpha_{i}^{-1}\beta_{i}^{-1}w_{2},$
$u=v=u^{-1}\beta_{i}u$\\\hline
\end{tabular}
$%
\end{tabular}
\]
where
\[
w_{2}=\prod_{k=i+1}^{\rho}[\alpha_{k},\beta_{k}]\prod_{k=1}^{j-1}\gamma_{k}%
\]

When the group $G$ is an additively written abelian group the formulas are
much simpler since all commutators will disappear and elements are equal to
their conjugates. To avoid confusing the additive formulas with the
multiplicative formulas in the general case we will use uppercase letters for
a generating vector viz., $(A_{1},\ldots A_{\rho},B_{1},\ldots B_{\rho}%
,C_{1},\ldots C_{\rho}).$ In this case an arbitrary assignment of elements of
$G$ to the generators $\eta:\alpha_{i}\rightarrow A_{i},$ $\beta
_{i}\rightarrow B_{i}$ will define an element of $\mathrm{Hom}(\Gamma,G)$
since the commutation relation is trivially satisfied. Our transformations
above may be written
\begin{align}
\mathfrak{A}_{i}^{k}  &  :B_{i}\rightarrow B_{i}+kA_{i}\label{eq-MCGgens-AB}\\
\mathfrak{B}_{i}^{k}  &  :A_{i}\rightarrow A_{i}+kB_{i}\nonumber\\
\mathfrak{Z}_{i}^{k}  &  :A_{i}\longrightarrow A_{i}+kA_{i+1},B_{i+1}%
\longrightarrow B_{i+1}-kB_{i}\nonumber\\
\mathfrak{R}_{i}  &  :A_{i}\rightarrow B_{i},\text{ }B_{i}\rightarrow
-A_{i}\nonumber\\
\mathfrak{S}_{i.}  &  :A_{i}\rightarrow A_{i+1},\text{ }B_{i}\rightarrow
B_{i+1},\text{ }A_{i+1}\rightarrow A_{i},\text{ }B_{i+1}\rightarrow
B_{i}\nonumber
\end{align}
and trivial on all other generators. If we include branch points then we must
also add an assignment $\eta:\gamma_{j}\rightarrow C_{j},$ with
\begin{equation}
o(C_{j})=m_{j}, \label{eq-GV-ab1}%
\end{equation}
and
\begin{equation}
C_{1}+C_{2}+\cdots+C_{r}=0. \label{eq-GV-ab2}%
\end{equation}
We can then add to our table of transformations
\begin{align}
\mathfrak{T}_{i}  &  :C_{j}\rightarrow C_{j+1},C_{j+1}\rightarrow C_{j},\text{
}m_{j}=m_{j+1}\label{eq-MCGgens-mixed}\\
\mathfrak{U}_{i,j}^{k}  &  :B_{i}\rightarrow B_{i}+kC_{j}\nonumber\\
\mathfrak{V}_{i,j}^{k}  &  :A_{i}\rightarrow A_{i}+kC_{j}\nonumber
\end{align}
Obviously, any permutation of the $C_{j}$ preserving order is permissible.

When $G$ is abelian then the induced action of $\mathrm{Aut}(\Gamma)$ on a
generating vector can be written in matrix vector form. Suppose that $\xi
^{-1}\in\mathrm{Aut}(\Gamma)$ is induced by a homeomorphism of $T^{o}$
preserving branch order. Then
\[
(A_{1},\ldots A_{\rho},B_{1},\ldots B_{\rho},C_{1},\ldots C_{r}%
)\longrightarrow(A_{1},\ldots A_{\rho},B_{1},\ldots B_{\rho},C_{1},\ldots
C_{r})M
\]
where $M$ is a matrix of the form
\begin{equation}
M=\left[
\begin{array}
[c]{cc}%
Sp & Z\\
0 & P
\end{array}
\right]  \label{eq-matdecomp}%
\end{equation}
where $Sp$ is a $2\rho\times2\rho$ integral symplectic matrix, $P$ is an
$r\times r$ permutation matrix, and $Z$ is arbitrary integer valued matrix.
The matrix $Sp$ is the induced symplectic automorphism on $H_{1}%
(T;\mathbb{Z}).$ Perhaps the easiest way to see this is that each of the
generators listed above has the given matrix decomposition.

\subsection{Hyperbolic - elliptic decomposition, abelian case}

In the group $\Gamma$ the generators $\left\langle \alpha_{i},\beta_{i}:1\leq
i\leq\rho\right\rangle $ are hyperbolic elements and the $\{\gamma_{j}:1\leq
j\leq r\}$ are elliptic elements. Accordingly we define the hyperbolic
(unramified) and elliptic (ramified) parts of $G:$ $G^{h}$ $=\left\langle
A_{i},B_{i}:1\leq i\leq\rho\right\rangle ,$ $G^{e}$ $=\left\langle C_{j}:1\leq
j\leq r\right\rangle .$ We are interested in the action of $\mathrm{Aut}%
(\Gamma)$ on the splitting of $G.$ According define for $\xi\in\mathrm{Aut}%
(\Gamma)$%
\[
G^{h,\xi}=\left\langle \eta\circ\xi^{-1}(\alpha_{i}),\eta\circ\xi^{-1}%
(\beta_{i}):1\leq i\leq\rho\right\rangle ,G^{e,\xi}=\left\langle \eta\circ
\xi^{-1}(\gamma_{j}):1\leq j\leq r\right\rangle .
\]
We have the following:

\begin{proposition}
Let $\eta:\alpha_{i}\rightarrow A_{i},$ $\beta_{i}\rightarrow B_{i},$
$\gamma_{j}\rightarrow C_{j},$ be an arbitrary assignment. Then $\eta$ extends
to an epimorphism $\eta:\Gamma\rightarrow G$ if and only if the relations
(\ref{eq-GV-ab1}) and (\ref{eq-GV-ab2}) hold and $G=\left\langle G^{h}%
,G^{e}\right\rangle .$
\end{proposition}

\begin{proposition}
Let $\eta:\Gamma\rightarrow G$ be an epimorphism then $G^{e}=G^{e,\xi}$ for
all $\xi\in\mathrm{Aut}(\Gamma).$
\end{proposition}

\begin{proof}
Let $\Gamma^{e}=\left\langle \gamma_{j}:1\leq j\leq r\right\rangle .$ If
$N(\Gamma^{e})$ is the normal closure of $\Gamma^{e}$ then $G^{e}=\eta\left(
N(\Gamma^{e})\right)  .$ It is well known that in $\Gamma$ every elliptic
element is conjugate to one of the $\gamma_{j},$ and hence $\left\langle
\eta\circ\xi^{-1}(C_{j}):1\leq j\leq r\right\rangle =$ $\eta\left(  \xi
^{-1}(N(\Gamma^{e}))\right)  =$ $\eta\left(  N(\Gamma^{e})\right)  =G^{e}%
.$\newline
\end{proof}

\begin{proposition}
The elliptic components of two abelian generating vectors ~$C_{j}:1\leq j\leq
r,$ $C_{j}^{\prime}:1\leq j\leq r,$ are $\mathcal{M}(\Gamma)$-equivalent if
and only if one is a permutation of the other: $C_{j}^{\prime}=C_{\pi j}$ for
some permutation $\pi.$
\end{proposition}

\begin{proof}
The automorphism $\xi$ is induced by a homeomorphism $h$ of $(T^{o},Q_{0})$
which must permute the punctures $Q_{1},\ldots Q_{r}.$ Since $h$ is
orientation-preserving then $h^{-1}(\widetilde{\gamma_{j}})$ is a conjugate of
some $\widetilde{\gamma_{\pi j}}$ for some permutation $\pi.$
\end{proof}

Next let $G^{e}$ be a subgroup of $G$ \ and consider all possible subgroups
$H$ such that $G=\left\langle H,G^{e}\right\rangle .$ Let $H_{0}$ be one of
smallest order. Then we may also assume that $G^{h}=H_{0}.$

\begin{proposition}
\label{pr-decomp}Let $\eta:\Gamma\rightarrow G$ be an action and let $H_{0}$
be the subgroup of smallest order such that $G=\left\langle H_{0}%
,G^{e}\right\rangle .$ Then there is $\xi$ $\in\mathcal{M}(\Gamma)$ such that
$G^{h,\xi}=H_{0}$ and $G^{e,\xi}=G^{e}.$
\end{proposition}

\begin{proof}
The equality $G^{e,\xi}=G^{e}$ is automatic. Each $A_{i}$ and $B_{i}$ may be
written
\begin{align*}
A_{i}  &  =A_{i}^{^{\prime}}+a_{i,1}C_{1}+\cdots+a_{i,r}C_{r},\text{ }%
A_{i}^{^{\prime}}\in H_{0}\\
B_{i}  &  =B_{i}^{^{\prime}}+b_{i,1}C_{1}+\cdots+b_{i,r}C_{r},\text{ }%
B_{i}^{^{\prime}}\in H_{0}%
\end{align*}
By applying $\mathfrak{V}_{i,j}^{-a_{i,j}}$ we can remove the component
$a_{i,j}C_{j}$ from the expression for $A_{i}$ and affect no other terms. Thus
all the $A_{i}$ can be reduced to $A_{i}^{^{\prime}}.$ Similarly the $B_{i}$
can be reduced to $B_{i}^{^{\prime}}.\ $Now $\left\langle A_{i}^{\prime}%
,B_{i}^{\prime}:1\leq i\leq\rho\right\rangle \subseteq H_{0}$ and
$G=\left\langle \left\langle A_{i}^{\prime},B_{i}^{\prime}:1\leq i\leq
\rho\right\rangle ,G^{e}\right\rangle $ by construction, so $H_{0}%
=\left\langle A_{i}^{\prime},B_{i}^{\prime}:1\leq i\leq\rho\right\rangle $ by minimality.
\end{proof}

\subsection{The elementary abelian case}

Finally let us assume that $G=\mathbb{F}_{p}^{w},$ an elementary abelian group
considered as a vector space over $\mathbb{F}_{p}.$ The generators
$A_{i},B_{i},C_{j}$ may then be thought of as vectors over $\mathbb{F}_{p}.$ Now let $X_{AB}$ be the matrix $\left[
\begin{array}
[c]{cccccc}%
A_{1} & \cdots & A_{\rho} & B_{1} & \cdots & B_{\rho}%
\end{array}
\right]  $ and $X_{C}$ be the matrix $\left[
\begin{array}
[c]{ccc}%
C_{1} & \cdots & C_{r}%
\end{array}
\right]  .$ Then by equation (\ref{eq-matdecomp}) the action of $(\omega,\xi)$
on $X=\left[
\begin{array}
[c]{cc}%
X_{AB} & X_{C}%
\end{array}
\right]  $ is:
\[
\left[
\begin{array}
[c]{cc}%
X_{AB} & X_{C}%
\end{array}
\right]  \rightarrow M_{\omega}\left[
\begin{array}
[c]{cc}%
X_{AB} & X_{C}%
\end{array}
\right]  M_{\xi}^{-1}%
\]
where $M_{\omega}\in GL_{w}(p)$ and
\[
M_{\xi}^{-1}=\left[
\begin{array}
[c]{cc}%
Sp & Z\\
0 & P
\end{array}
\right]
\]
Next by Proposition \ref{pr-decomp} we assume that $\eta$ is chosen so that
$\mathbb{F}_{p}^{w}=G^{h}\oplus G^{e}.$ Now write $\mathbb{F}_{p}%
^{w}=\mathbb{F}_{p}^{u}\oplus\mathbb{F}_{p}^{v}$ and find an $N$ $\in
GL_{w}(p)$ such $NG^{h}=\mathbb{F}_{p}^{u}\oplus0$ and $NG^{e}=0\oplus
\mathbb{F}_{p}^{v}.$ It then follows that
\begin{equation}
N\left[
\begin{array}
[c]{cc}%
X_{AB} & X_{C}%
\end{array}
\right]  =\left[
\begin{array}
[c]{cc}%
Y_{AB} & 0\\
0 & Y_{C}%
\end{array}
\right]  \label{eq-canonicalGV}%
\end{equation}
where $Y_{AB}$ is $u\times2\rho$ matrix of rank $u$ and $Y_{C}$ is a $v\times
r$ matrix of rank $v.$ Thus we may assume the generating vector $X$ of our
initial $\eta$ has the form of the right hand side of equation
(\ref{eq-canonicalGV}). Now suppose that we look at a transformed generating
vector $M_{\omega}XM_{\xi}^{-1}$ that has the same form i.e.,%
\begin{align*}
\left[
\begin{array}
[c]{cc}%
W_{AB} & 0\\
0 & W_{C}%
\end{array}
\right]   &  =\left[
\begin{array}
[c]{cc}%
M_{11} & M_{12}\\
M_{21} & M_{22}%
\end{array}
\right]  \left[
\begin{array}
[c]{cc}%
Y_{AB} & 0\\
0 & Y_{C}%
\end{array}
\right]  \left[
\begin{array}
[c]{cc}%
Sp & Z\\
0 & P
\end{array}
\right] \\
&  =\left[
\begin{array}
[c]{cc}%
M_{11}Y_{AB}Sp & M_{11}Y_{AB}Z+M_{12}Y_{C}P\\
M_{21}Y_{AB}Sp & M_{21}Y_{AB}Z+M_{22}Y_{C}P
\end{array}
\right]
\end{align*}
Since $M_{21}Y_{AB}Sp=0$ and $Y_{AB}Sp$ has rank $u$ then $M_{21}=0,$and
$M_{11}$ and $M_{22}$ are invertible. Since $M_{11}Y_{AB}Z+M_{12}Y_{C}P$ is
also assumed zero, it follows that
\begin{align*}
\left[
\begin{array}
[c]{cc}%
W_{AB} & 0\\
0 & W_{C}%
\end{array}
\right]   &  =\left[
\begin{array}
[c]{cc}%
M_{11}Y_{AB}Sp & 0\\
0 & M_{22}Y_{C}P
\end{array}
\right] \\
&  =\left[
\begin{array}
[c]{cc}%
M_{11} & \\
& M_{22}%
\end{array}
\right]  \left[
\begin{array}
[c]{cc}%
Y_{AB} & 0\\
0 & Y_{C}%
\end{array}
\right]  \left[
\begin{array}
[c]{cc}%
Sp & \\
0 & P
\end{array}
\right]
\end{align*}
i.e., we can just find the equivalence classes of the unramified and the
ramified components independently. We summarize the results in a theorem.

\begin{theorem}
\label{th-split}Suppose that the group $G=\mathbb{F}_{p}^{w}$ acts on a
surface $S$ with signature $(\rho:p^{r}),$ and genus $1+p^{w}(\rho
-1)+p^{w-1}\frac{r\left(  p-1\right)  }{2}$ determined by an epimorphism
$\eta:\Gamma\rightarrow G$ (note that $r=1$ is infeasible). Then there are unique
integers $u,v$ such that $w=u+v,$ $0\leq u\leq2\rho,$ $1\leq v<r,$ a
$u\times2\rho$ matrix $W_{AB}$ of rank $u,$ and a $v\times r$ matrix $W_{C}$
of rank $v$ such that the matrix $\left[
\begin{array}
[c]{cc}%
W_{AB} & 0\\
0 & W_{C}%
\end{array}
\right]  $ is a representative of the topological equivalence class of
$G$-actions on $S.$ Moreover $\left[
\begin{array}
[c]{cc}%
W_{AB}^{\prime} & 0\\
0 & W_{C}^{\prime}%
\end{array}
\right]  $ is another representative satisfying similar rank conditions for
$W_{AB}^{\prime}$ and $W_{C}^{\prime}$ then $W_{AB}$ and $W_{AB}^{\prime}$
define equivalent unramified actions with signature $(\rho:-)$ and $W_{C}$ and
$W_{C}^{\prime}$ define equivalent purely ramified actions with signature
$(0:p^{r}).$
\end{theorem}

\begin{corollary}
\label{cor-counting}For $0\leq u\leq2\rho,$ let $h_{u}$ be the number of
equivalence classes of unramified actions of $\mathbb{F}_{p}^{u}$ on a surface
of genus $1+p^{u}(\rho-1)$ with signature $(\rho:-)$ and $h_{u}=0$ otherwise.
For $1\leq v<r$ let $e_{v}$ be the number of equivalence classes of purely
ramified actions of $\mathbb{F}_{p}^{v}$ on a surface $S$ of genus
$1+p^{w-1}\frac{r\left(  p-1\right)  }{2}-p^{w}$ with signature $(0:p^{r})$
and set $e_{v}=0$ otherwise. Then the number of inequivalent actions of
$\mathbb{F}_{p}^{w}$ on a surface of genus $1+p^{w}(\rho-1)+p^{w-1}%
\frac{r\left(  p-1\right)  }{2}$ with signature $(\rho:p^{r})$ is given by%
\[
\text{\#\textrm{actions} }=\sum\limits_{u=0}^{n}h_{u}e_{w-u}.
\]

\end{corollary}

\subsection{Cohomological invariants\label{subsec-cohominv}}

For abelian covers we may use cohomology to concoct an $\mathcal{M}$-invariant
to distinguish classes. As noted above, for an abelian group $G$,
$\mathrm{Hom}(\Gamma,G)$ classifies covers $T$ by subgroups of $G.$ We have
the following sequence of equivalences%
\[
\mathrm{Hom}(\Gamma,G)\backsimeq\mathrm{Hom}(H_{1}(T;\mathbb{Z}),G)\backsimeq
H^{1}(T;G).
\]
Given two elements of $\eta_{1},\eta_{2}\in H^{1}(T;G)$ we may consider the
cup product $\eta_{1}\cup\eta_{2}\in H^{1}(T;G\otimes G).$ If $h$ is an
orientation-preserving homeomorphism then
\[
h^{\ast}\eta_{1}\cup h^{\ast}\eta_{2}=h^{\ast}(\eta_{1}\cup\eta_{2})=\eta
_{1}\cup\eta_{2}.
\]
In particular if $\eta_{1}=\eta_{2}=\eta$ then%
\[
h^{\ast}\eta\cup h^{\ast}\eta=\eta\cup\eta
\]
since $h$ is orientation preserving. It is not hard to show that if
$\eta:\alpha_{i}\rightarrow A_{i},$ $\eta:\beta_{i}\rightarrow B_{i},$ then
\begin{equation}
\eta\cup\eta=\sum\limits_{i=1}^{\rho}\left(  A_{i}\otimes B_{i}-B_{i}\otimes
A_{i}\right)  \label{eq-CupInv}%
\end{equation}
Observe that the invariance of $\eta\cup\eta$ under the transforms in
Subsection \ref{subsec-AGxAGam} may also be proven by direct computation.

\begin{example}
\label{ex-cup}Let $G$ = $C_{n}\oplus C_{m}$  with $x$
generating $C_{n}$ and $y$ generating $C_{m}$ additively (cyclic groups of order $m$ and $n$).
Consider the following epimorphisms%
\begin{align*}
\eta_{1}  &  :\left(  \alpha_{1},\beta_{1},\alpha_{2},\beta_{2}\right)
\rightarrow\left(  x,y,0,0\right) \\
\eta_{2}  &  :\left(  \alpha_{1},\beta_{1},\alpha_{2},\beta_{2}\right)
\rightarrow\left(  x,0,y,0\right)
\end{align*}
Then
\begin{align*}
\eta_{1}\cup\eta_{1}  &  =x\otimes y-y\otimes x\\
\eta_{2}\cup\eta_{2}  &  =0
\end{align*}
and so the epimorphisms are different under the action of \textrm{Aut}%
$(\Gamma(2;-))$.
\end{example}

\section{\label{Sec-Unramified}The unramified abelian case}

In this section, we consider the case when $G$ is abelian and $\Gamma
=\Pi_{\rho}$ is torsion free. For the special case when $G$ is elementary
abelian, complete results were derived in \cite{Bro4}. We shall apply the
ideas from this case to the more general abelian case. Due to the complexity
of the problem, we shall only produce partial results, though the results we
derive will be sufficient to produce explicit results for genus up to $65$. We
shall also illustrate how in principle one could use the results to classify
all fixed point abelian actions for arbitrary genus. It should be noted that
we do not consider the abelian case for general $\Gamma$ because Theorem
\ref{th-split} no longer holds so there could be overlap between $G^{h}$ and
$G^{e}$ making the problem much more difficult.

\subsection{\label{subsec-unramified}The unramified elementary abelian case}

The following was proved in \cite{Bro4} and completely classifies all
unramified elementary abelian actions up to topological equivalence (by $p$-rank of an abelian group $G$, we mean the number of invariant factors of $G$). The
corollary immediately follows.

\begin{thm}
\label{thm-unramified} Suppose $\Gamma=\Pi_{\rho}$ is a surface group of orbit
genus $\rho$ and generators $\alpha_{1},\dots,\alpha_{\rho},\beta_{1},\dots
b_{\rho}$ and $G$ is an elementary abelian group of $p$-rank $r\leq2\rho$ with
generators $\omega_{1},\dots,\omega_{r}$. Then there exists an integer
$r/2\leq K\leq\mathrm{{min}(\rho,r)}$ such that any epimorphism from $\Gamma$
onto $G$ is $\aut{(G)}\times\aut{(\Gamma )}$-equivalent to one of those below:%

\[
\eta:=%
\begin{cases}
\eta(\alpha_{i})=\omega_{i} & \text{$i\leq K$}\\
\eta(\alpha_{i})=0 & \text{$i>K$}\\
\eta(\beta_{i})=\omega_{i+K} & \text{$i\leq r-K$}\\
\eta(\beta_{i})=0 & \text{$i>r-K$}%
\end{cases}
\]

\end{thm}

\begin{corollary}
\label{cor-Number} Suppose $G$ is an elementary abelian group of $p$-rank $r$, let
$\mathcal{M}_{\sigma}$ denote the mapping class group of a closed surface of
genus $\sigma$ and let $\rho=(\sigma-1+p^{r})/p^{r}$. If $\rho$ is not an
integer, there are no conjugacy classes of subgroups of $\mathcal{M}_{\sigma}$
isomorphic to $G$ with fixed point free action. Else, the number of conjugacy
classes of subgroups of $\mathcal{M}_{\sigma}$ isomorphic to $G$ with fixed
point free action is calculated as follows: if $r\leq\rho$, there are $r/2+1$
different classes of epimorphisms for $r$ even and $(r+1)/2$ epimorphisms if
$r$ is odd, if $r=\rho+i$ with $0<i<\rho$ there are $(\rho-i)/2$ classes of
epimorphisms if $\rho-i$ is even and $(\rho-i+1)/2$ if $\rho-i$ is odd, and if
$r=2\rho$ or $r=1$, there is just one class.
\end{corollary}

We illustrate with an example.

\begin{example}
Suppose $\Gamma$ has genus $2$ and $G$ is an elementary abelian $p$-group of
$p$-rank $2$. If $\eta\colon\Gamma\rightarrow G$ is a surface kernel epimorphism,
then using the Riemann-Hurwitz formula, the genus $\sigma$ of the kernel will
be $\sigma=p^{2}+1$. In fact, by simple application of the Riemann Hurwitz
formula, it can be shown that for $p>5$, $\Gamma$ is the only Fuchsian group
which admits an elementary abelian quotient of order $p^{2}$ with kernel of
orbit genus $p^{2}+1$ (for $p=5$, see Example \ref{ex-genus26}). This means
for $p>5$, in $\mathcal{M}_{p^{2}+1}$, the conjugacy classes of elementary
abelian subgroups of order $p^{2}$ will be in $1-1$ correspondence with the
classes of epimorphisms from $\Gamma$ onto $G$. Applying Corollary
\ref{cor-Number}, there are two such classes.
\end{example}

\subsection{\label{subsec-genus2}The general unramified abelian case}

We now consider partial results for the general abelian case. Before we start,
we introduce some notation and terminology.

\begin{notation}
For the rest of this section, $G$ will denote an additively written abelian
group of $p$-rank $r$ with invariant factors $n_{1}\geq n_{2}\geq\dots\geq n_{r}$
where $n_{i+1}|n_{i}$ and $\omega_{1},\dots,\omega_{r}$ are a fixed set of
generators of orders $n_{1},\dots,n_{r}$ respectively. Also, for an integer
$n$, $C_{n}$ denotes the cyclic group of order $n$.
\end{notation}

The following two steps can be taken to classify all fixed point free
$G$-actions on a surface - first determine a set of epimorphisms with the
property that every epimorphism is equivalent to one in this set, and then
reduce this set so no two epimorphisms are equivalent. In the following, we
shall consider the first step. Though in general we shall not consider the
problem of distinguishing between classes, we shall present explicit examples
showing how in principle one could tackle this problem. Note that $r\leq2\rho
$, so we only need consider epimorphisms from $\Gamma$ to $G$ with $p$-rank at
most $2\rho$.

\begin{proposition}
\label{prop-abclasses} Any epimorphism $\eta\colon\Gamma\rightarrow G$ is
equivalent to one of those described below:

\begin{enumerate}
\item $r\leq\rho$.%

\[
\eta:=%
\begin{cases}
\eta(\alpha_{i})=\omega_{i} & i\leq r\\
\eta(\alpha_{i})=0 & i>r\\
\eta(\beta_{i})=\sum\limits_{j=i+1}^{r}N_{\beta_{i},\omega_{j}}\omega_{j} &
i\leq r\\
\eta(\beta_{i})=0 & i>r
\end{cases}
\]
where

\begin{itemize}
\item $(n_{i+1},N_{\beta_{i},\omega_{i+1}})>1$ or $N_{\beta_{i},\omega_{i+1}%
}=1$ for each $i\leq r-1$

\item for each $i$, if $j>i+1$ is the largest integer such that $N_{\beta
_{i},\omega_{j}}\neq0$, then $\langle N_{\beta_{i},\omega_{j}}\omega
_{j}\rangle\leq\langle N_{\beta_{i+1},\omega_{j}}\omega_{j}\rangle$.
\end{itemize}

\item $r>\rho$.%

\[
\eta:=%
\begin{cases}
\eta(\alpha_{i})=\omega_{i} & 1\leq i\leq\rho\\
\eta(\beta_{i})=\sum\limits_{j=i+1}^{r}N_{\beta_{i},\omega_{j}}\omega_{j}, &
1\leq i\leq2\rho-r\\
\eta(\beta_{i})=\omega_{2\rho-i+1}+\sum\limits_{j=i+1}^{2\rho-i}N_{\beta
_{i},\omega_{j}}\omega_{j} & 2\rho-r<i<\rho\\
\eta(\beta_{\rho})=\omega_{\rho+1} & i=\rho
\end{cases}
\]
where

\begin{itemize}
\item for $i>2\rho-r$ either the order of $N_{\beta_{i},\omega_{j}}\omega_{j}$
does not divide $n_{r}$ or $N_{\beta_{i},\omega_{j}}=0$

\item for $i\leq2\rho-r$ if $j>i+1$ is the largest integer such that
$N_{\beta_{i},\omega_{j}}\neq0$, then $\langle N_{\beta_{i},\omega_{j}}%
\omega_{j}\rangle\leq\langle N_{\beta_{i+1},\omega_{j}}\omega_{j}\rangle$

\item for $i\leq2\rho-r-1$, $(n_{i},N_{\beta_{i},\omega_{i+1}})>1$ or
$N_{\beta_{i},\omega_{i+1}}=1$.
\end{itemize}
\end{enumerate}
\end{proposition}

\begin{proof}
We shall use induction on the $p$-rank of the group $G$. For $p$-rank $r=1$, we are
done by Theorem 14 of \cite{Har2}. Assuming the result holds for $p$-rank $r-1$,
we shall prove it holds for $p$-rank $r$. The proof falls into two different cases
depending upon whether $1<r\leq\rho$ or $\rho<r\leq2\rho$. First suppose that
$r\leq\rho$.

Let $\eta\colon\Gamma\rightarrow G=C_{n_{1}}\times\dots\times C_{n_{r}}$
denote the epimorphism onto $G$ and $\Phi\colon G\rightarrow C_{n_{1}}%
\times\dots\times C_{n_{r-1}}$ the projection map onto the quotient group
$G/C_{n_{r}}$ (by abuse of notation, we identify $C_{n_{1}}\times\dots\times
C_{n_{r-1}}$ with the corresponding subgroup of $G$). Since $\Gamma$ is
torsion free, all subgroups will be torsion free. In particular, the map
$\Phi\circ\eta$ will be a surface kernel epimorphism from $\Gamma$ onto
$C_{n_{1}}\times\dots\times C_{n_{r-1}}$, so by induction, will be equivalent
to one as given in the statement of the proposition. Lifting to $G$ and
composing with appropriate automorphisms of $G$, it follows that $\eta$ is
equivalent to an epimorphism of the following form:%

\[
\eta:=%
\begin{cases}
\eta(\alpha_{i})=\omega_{i} & \text{$i<r$}\\
\eta(\alpha_{i})=m_{i}\omega_{r} & \text{$i\geq r$}\\
\eta(\beta_{i})=M_{i}\omega_{r}+\sum\limits_{j=i+1}^{r-1}N_{\beta_{i}%
,\omega_{j}}\omega_{j} & \text{$i\leq r$}\\
\eta(\beta_{i})=M_{i}\omega_{r} & \text{$i>r$}%
\end{cases}
\]

\noindent for integers $m_{i}$, $r\leq i\leq\rho$ and $M_{j}$, $1\leq
j\leq\rho$. Similar to the elementary abelian case in \cite{Bro4}, we now
reduce using the automorphisms from $\aut{(\Gamma)}$ developed in Section
\ref{subsec-AGxAGam} and the automorphisms of the abelian group $G$. In most
instances, since the reduction steps are between pairs of pairs of generators
and constant on all others, we shall just consider the pairs which are
changed. To start, let $k$ be the smallest integer such that $(M_{2}%
-kM_{1})\omega_{r}$ generates the subgroup generated by $M_{1}\omega_{r}$ and
$M_{2}\omega_{r}$ and let $\Phi\in\aut(G)$ where $\Phi(\omega_{1})=\omega
_{1}+k\omega_{2}$ and the identity map on all other generators. Then the map
$\mathfrak{A}_{1}^{n_{2}-N_{\beta_{1},\omega_{2}}}\Phi\circ\mathfrak{Z}%
_{1}^{k}$ modifies the images of $\alpha_{1}$, $\alpha_{2}$, $\beta_{1}$ and
$\beta_{2}$ but acts trivially on the images of $\alpha_{i}$ and $\beta_{i}$
for $i\geq3$. Specifically, after renaming the coefficients, we get:
\begin{align*}
(\alpha_{1},\beta_{1},\alpha_{2},\beta_{2})  &  \longrightarrow(\omega
_{1},M_{1}\omega_{r}+\sum\limits_{j=2}^{r-1}N_{\beta_{i},\omega_{j}}\omega
_{j},\omega_{2},M_{2}\omega_{r}+\sum\limits_{j=3}^{r-1}N_{\beta_{i},\omega
_{j}}\omega_{j})\\
&  \longrightarrow(\omega_{1},\sum\limits_{j=2}^{r}N_{\beta_{1},\omega_{i}%
}\omega_{i},\omega_{2},\sum\limits_{j=3}^{r}N_{\beta_{2},\omega_{i}}\omega
_{i})
\end{align*}
where $\langle N_{\beta_{1},r}\omega_{r}\rangle\leq\langle N_{\beta_{2}%
,r}\omega_{r}\rangle$. Applying similar transformations to all proceeding
pairs of pairs up to the pair $(\alpha_{r},\beta_{r})$ we get:%

\[
\eta:=%
\begin{cases}
\eta(\alpha_{i})=\omega_{i} & \text{$i<r$}\\
\eta(\alpha_{i})=m_{i}\omega_{r} & \text{$i\geq r$}\\
\eta(\beta_{i})=\sum\limits_{j=i+1}^{r}N_{\beta_{i},\omega_{j}}\omega_{j} &
\text{$i\leq r$}\\
\eta(\beta_{i})=M_{i}\omega_{r} & \text{$i>r$}%
\end{cases}
\]
where $\langle N_{\beta_{i},r}\omega_{r}\rangle\leq\langle N_{\beta_{i+1}%
,r}\omega_{r}\rangle$ for all $i$. In particular, this implies
\[
\langle N_{\beta_{1},r},N_{\beta_{2},r},\dots,N_{\beta_{r-1},r}\omega
_{r}\rangle\leq\langle N_{\beta_{r},r}\omega_{r}\rangle.
\]

Next, we reduce the pairs $(\alpha_{i},\beta_{i})$ with $i\geq r$. First, to
any such pair, we apply $\mathfrak{B}_{i}^{k}$ where $k$ is the smallest
integer such that $(m_{i}+kM_{i})\omega_{r}$ generates the subgroup generated
by $m_{i}\omega_{r}$ and $M_{i}\omega_{r}$ (for brevity, by abuse of notation,
we rename $(m_{i}+kM_{i})$ by $m_{i}$). Then for a given $i$, since $\langle
M_{i}\omega_{r}\rangle\leq\langle m_{i}\omega_{r}\rangle$, there exists $n$
such that $(M_{i}+nm_{i})\omega_{r}=0$. Applying the automorphism
$\mathfrak{A}_{i}^{n}$, we get $(\alpha_{i},\beta_{i})\rightarrow
(0,M_{i}\omega_{r})$. Following this, for the pairs, $(\alpha_{\rho-1}%
,\beta_{\rho-1},\alpha_{\rho},\beta_{\rho})$, we apply $\mathfrak{A}_{\rho
-1}^{k}$ where $k$ is the smallest integer such that $(M_{\rho}-kM_{\rho
-1})\omega_{r}$ generates the subgroup generated by $M_{\rho}\omega_{r}$ and
$M_{\rho-1}\omega_{r}$. Renaming $(M_{\rho}-kM_{\rho-1})\omega_{r}$ by
$M_{\rho}$, we get $(\alpha_{\rho-1},\beta_{\rho-1},\alpha_{\rho},\beta_{\rho
})\rightarrow(M_{\rho-1}\omega_{r},0,M_{\rho}\omega_{r},0)$ where $\langle
M_{\rho-1}\omega_{r}\rangle\leq\langle M_{\rho}\omega_{r}\rangle$. Finally, we
eliminate the image of $\beta_{\rho-1}$ by applying $\mathfrak{A}_{\rho-1}%
^{n}$ where $n$ is the smallest integer such that $(M_{\rho-1}-nM_{\rho
})\omega_{r}=0$ and then apply $\mathfrak{S}_{\rho-1}$ and $\mathfrak{R}%
_{\rho-1}$ giving $(\alpha_{\rho-1},\beta_{\rho-1},\alpha_{\rho},\beta_{\rho
})\rightarrow(0,0,M_{\rho-1}\omega_{r},0)$ (after renaming coefficients).
Applying similar transformations to pairs $(\alpha_{i},\beta_{i},\alpha
_{i+1},\beta_{i+1})$ with $i>r$, we get:%

\[
\eta:=%
\begin{cases}
\eta(\alpha_{i})=\omega_{i} & \text{$i<r$}\\
\eta(\alpha_{r})=M\omega_{r} & \text{$i=r$}\\
\eta(\alpha_{i})=0 & \text{$i>r$}\\
\eta(\beta_{i})=\sum\limits_{j=i+1}^{r}N_{\beta_{i},\omega_{j}}\omega_{j} &
\text{$i\leq r$}\\
\eta(\beta_{i})=0 & \text{$i>r$}%
\end{cases}
\]
for some integer $M$. Finally, we apply $\mathfrak{Z}_{r-1}^{k}$ where $k$ is
the smallest integer such that $(M-kN_{\beta_{r-1},\omega_{r}},n_{r})=1$ (note
that such an integer exists since $M\omega_{r}$ and $N_{\beta_{r-1},\omega
_{r}}\omega_{r}$ generate $\langle\omega_{r}\rangle$). Finally, by applying
$\mathfrak{G}_{r}$ and an appropriate automorphism of $G$, after renaming the
coefficients, we get:%

\[
\eta:=%
\begin{cases}
\eta(\alpha_{i})=\omega_{i} & \text{$i\leq r$}\\
\eta(\alpha_{i})=0 & \text{$i>r$}\\
\eta(\beta_{i})=\sum\limits_{j=i+1}^{r}N_{\beta_{i},\omega_{j}}\omega_{j} &
\text{$i\leq r$}\\
\eta(\beta_{i})=0 & \text{$i>r$}%
\end{cases}
\]

For the last reduction step, if $(N_{\beta_{r-1} ,\omega_{r}} ,n_{r})=1$, we
first apply $\Phi\in\aut (G)$ such that $\Phi(N_{\beta_{r-1} ,\omega_{r}}
\omega_{r} )=\omega_{r}$ and identity on all other generators. Then, assuming
$b$ is the integer with $\Phi(\omega_{r} )=b\omega_{r}$, we perform
$\mathfrak{G}_{r} \circ\mathfrak{Z}_{r-1}^{b+1} \circ\mathfrak{G}_{r}
\circ\Phi$ giving $(\alpha_{r-1} ,\beta_{r-1} ,\alpha_{r} ,\beta_{r}
)\rightarrow(\omega_{r-1} ,\omega_{r} ,\omega_{r} ,0)$. Note that for $l<r$,
none of the transformations used during this proof change the fact that
$(n_{l}, N_{\beta_{l-1},\omega_{l}})>1$ or $N_{\beta_{l-1},\omega_{l}}=1$.
Observe also that the reduction methods used imply that if $j>i+1$ is the
largest integer such that $N_{\beta_{i},\omega_{j}} \neq0$, then $\langle
N_{\beta_{i},\omega_{j}} \omega_{j} \rangle\leq\langle N_{\beta_{i+1}%
,\omega_{j}} \omega_{j} \rangle$.

Now suppose that $r>\rho$. Since the arguments regarding transformations are
similar to those used for the previous case, we shall skip many steps.
Assuming the result holds when the $p$-rank of the group is $r-1\geq\rho$,
induction implies after composing with appropriate automorphisms of $G$,
$\eta$ has the following form:%

\[
\eta:=%
\begin{cases}
\eta(\alpha_{i})=\omega_{i} & 1\leq i\leq\rho\\
\eta(\beta_{i})=M_{i}\omega_{r}+\sum\limits_{j=i+1}^{r-1}N_{\beta_{i}%
,\omega_{j}}\omega_{j}, & 1\leq i\leq2\rho-r+1\\
\eta(\beta_{i})=\omega_{2\rho-i+1}+\sum\limits_{j=i+1}^{2\rho-i}N_{\beta
_{i},\omega_{j}}\omega_{j} & 2\rho-r+1<i<\rho\\
\eta(\beta_{\rho})=\omega_{\rho+1} &
\end{cases}
\]

As with the previous case, for each $i<r$, we can use the transformations
$\mathfrak{Z}_{i}$, $\mathfrak{A}_{i}$ and appropriate automorphisms of $G$ to
move a generator for $\langle\omega_{r}\rangle$ to the image of $\beta
_{2\rho-r+1}$ giving $\eta(\beta_{2\rho-r+1})=\omega_{2\rho-r+1}%
+\sum\limits_{j=r+1}^{2\rho-r}N_{\beta_{i},\omega_{j}}\omega_{j}$. Applying
these transformations, we get%

\[
\eta:=%
\begin{cases}
\eta(\alpha_{i})=\omega_{i} & 1\leq i\leq\rho\\
\eta(\beta_{i})=\sum\limits_{j=i+1}^{r}N_{\beta_{i},\omega_{j}}\omega_{j}, &
1\leq i\leq2\rho-r\\
\eta(\beta_{i})=\omega_{2\rho-i+1}+\sum\limits_{j=i+1}^{2\rho-i}N_{\beta
_{i},\omega_{j}}\omega_{j} & 2\rho-r+1\leq i<\rho\\
\eta(\beta_{\rho})=\omega_{\rho+1} &
\end{cases}
\]
where $\langle N_{\beta_{1},\omega_{r}}\omega_{r}\rangle\leq\langle
N_{\beta_{2},\omega_{j}}\omega_{j}\rangle\leq\dots\leq\langle N_{\beta
_{2\rho-r+1},\omega_{r}}\omega_{r}\rangle$.

For the last reduction step, we apply an appropriate automorphism of $G$ that
acts trivially on all generators except $\omega_{r}$ and eliminates all
elements from the sum $\sum\limits_{j=i+1}^{r-1} N_{\beta_{i},\omega_{j}}
\omega_{j}$ with the order of $N_{\beta_{i},\omega_{j}} \omega_{j}$ dividing
$n_{r}$. Then we get:%

\[
\eta:=%
\begin{cases}
\eta(\alpha_{i})=\omega_{i} & 1\leq i\leq\rho\\
\eta(\beta_{i})=\sum\limits_{j=i+1}^{r}N_{\beta_{i},\omega_{j}}\omega_{j}, &
1\leq i\leq2\rho-r\\
\eta(\beta_{i})=\omega_{2\rho-i+1}+\sum\limits_{j=i+1}^{2\rho-i}N_{\beta
_{i},\omega_{j}}\omega_{j} & 2\rho-r<i
\end{cases}
\]
where the order of $N_{\beta_{i},\omega_{j}}\omega_{j}$ does not divide
$n_{r}$ for $i>2\rho-r$. As with the last case, if $i\leq2\rho-r$, the
reduction methods used implies that if $j>i+1$ is the largest integer such
that $N_{\beta_{i},\omega_{j}}\neq0$, then $\langle N_{\beta_{i},\omega_{j}%
}\omega_{j}\rangle\leq\langle N_{\beta_{i+1},\omega_{j}}\omega_{j}\rangle$.
\end{proof}

We note that though more refined sets of epimorphisms can be obtained (in
particular, the elementary abelian case), our goal was to provide a set of
epimorphisms independent of the genus and invariant factors. We finish with an
explicit example showing how to use Proposition \ref{prop-abclasses} to
determine all classes of fixed point free actions.


\begin{example}
\label{ex-genus33} Suppose that $\Gamma$ has genus $2$ and $G$ has $p$-rank $3$
with invariant factors $4,4,2$. Proposition \ref{prop-abclasses} implies there
are up to six classes of epimorphisms of the form $(\alpha_{1} ,\beta_{1}
,\alpha_{2} ,\beta_{2} )\rightarrow(\omega_{1} ,a\omega_{2} +b\omega_{3}
,\omega_{2} ,\omega_{3} )$ where $a=0,1,2$ and $b=0,1$. Using elements of
$\aut{(G)} \times\aut{(\Gamma )}$ it is fairly straight forward to reduce each
of these epimorphisms to one of the following three:%

\begin{tabular}
[c]{l}%
$\eta_{1}:=(\alpha_{1},\beta_{1},\alpha_{2},\beta_{2})\rightarrow(\omega
_{1},\omega_{2},\omega_{3},0)$\\
$\eta_{2}:=(\alpha_{1},\beta_{1},\alpha_{2},\beta_{2})\rightarrow(\omega
_{1},\omega_{3},\omega_{2},0)$\\
$\eta_{3}:=(\alpha_{1},\beta_{1},\alpha_{2},\beta_{2})\rightarrow(\omega
_{1},\omega_{2},2\omega_{2},\omega_{3})$%
\end{tabular}

\noindent We shall show that these three epimorphisms are distinct.

First note that if $\eta_{i}$ is equivalent to $\eta_{j}$, then there exists
$\theta\in\aut{(G)}$ and $\zeta\in\aut{(\Gamma )}$ such that $\theta\circ
\eta_{i}\circ\zeta=\eta_{j}$. Since $\cup$ is an $\aut{(\Gamma )}$ invariant,
it follows that $(\theta\circ\eta_{i})\cup(\theta\circ\eta_{i})=\eta_{j}%
\cup\eta_{j}$. Therefore, in order to show distinctness, it suffices to show
that for $i\neq j$, there does not exist $\theta\in\aut{(G)}$ such that
$(\theta\circ\eta_{i})\cup(\theta\circ\eta_{i})=\eta_{j}\cup\eta_{j}$. Since
$\aut{(G)}$ is finite, this can be checked directly. We illustrate by showing
that $\eta_{1}$ and $\eta_{2}$ are not equivalent.

First, taking the cup product we get
\[
\eta_{1}\cup\eta_{1}=\omega_{1}\otimes\omega_{2}-\omega_{2}\otimes\omega_{1}%
\]
and
\[
\theta\circ\eta_{2}\cup\theta\circ\eta_{2}=\theta(\omega_{1})\otimes
\theta(\omega_{3})-\theta(\omega_{3})\otimes\theta(\omega_{1}).
\]
Assuming that $\rho(\omega_{1})=A\omega_{1}+B\omega_{2}+C\omega_{3}$ and
$\rho(\omega_{3})=D\omega_{1}+E\omega_{2}+F\omega_{3}$, by simplifying the
tensor product and combining like terms, in order for these products to be
equal, we must have $(AE-DB)=1\mod{(4)}$. However, this is impossible since
both $D$ and $E$ must be divisible by $2$ (since $\omega_{3}$ has order $2$).
Thus $\eta_{1}$ and $\eta_{2}$ are not equivalent. All other cases follow
similar arguments and so it follows that there are three topological
equivalence classes of fixed point $G$-actions on a surface of genus $33$.
\end{example}

\subsection{\label{subsec-smallclass} Classification for Genera $\leq65$}

Due to the ad-hoc style arguments which seem required to distinguish classes,
we fall short of a general classification for abelian groups. However,
Proposition \ref{prop-abclasses} coupled with the following results will allow
a classification up to genus $65$ (and in fact is enough to classify certain
infinite families).

\begin{lemma}
\label{lem-equivclasses} If $G$ has $p$-rank $2$ and $n_{2}$ is square free, then
for each divisor $N$ of $n_{2}$ with $1\leq N\leq n_{2}$, any epimorphism
$\eta\colon\Gamma\rightarrow G$ is equivalent to one of those described below:%

\[
\eta:=%
\begin{cases}
\eta(\alpha_{1})=\omega_{1} & \\
\eta(\alpha_{2})=\omega_{2} & \\
\eta(\alpha_{i})=0 & i>2\\
\eta(\beta_{1})=N\omega_{2} & \\
\eta(\beta_{i})=0 & i>1
\end{cases}
\]

\noindent Moreover, no two such epimorphisms are equivalent under the action
of $\aut{(G)} \times\aut{(\Gamma )}$, so in particular, there are $d (n_{2})$
epimorphisms up to the action of $\aut{(G)} \times\aut{(\Gamma )}$ (where $d
(n_{2})$ denotes the number of divisors of $n_{2}$).
\end{lemma}

\begin{proof}
By Proposition \ref{prop-abclasses}, any epimorphism from $\Gamma$ onto a $p$-rank
$2$ abelian group will be equivalent to one of the form $\eta_{N}
\colon(\alpha_{1} ,\beta_{1} ,\alpha_{2} ,\beta_{2} )\rightarrow(\omega_{1}
,N\omega_{2} ,\omega_{2} ,0)$ with the image of all other generators being
trivial. Therefore, we just need to determine when $\eta_{N}$ and $\eta_{M}$
are equivalent for $M\neq N$. We start by showing that if $\gcd{(n,N)}
\neq\gcd{(n,M)}$, then $\eta_{N}$ and $\eta_{M}$ are not equivalent.

First, using the arguments from Example \ref{ex-genus33}, it suffices to show
that there does not exist $\theta\in\aut{(G)}$ such that $(\theta\circ\eta
_{N})\cup(\theta\circ\eta_{N})=\eta_{M}\cup\eta_{M}$. In general, if
$\theta\in\aut{(G)}$, then $\theta$ can be identified with a $2\times2$ matrix
with integer coefficients with the action of $\theta$ on the generators
$\omega_{1}$ and $\omega_{2}$ defined as follows:%

\[
\theta\cdot%
\begin{array}
[c]{|c|}%
\omega_{1}\\
\omega_{2}%
\end{array}
=%
\begin{array}
[c]{|cc|}%
A & B\\
C & D
\end{array}
\cdot%
\begin{array}
[c]{|c|}%
\omega_{1}\\
\omega_{2}%
\end{array}
=%
\begin{array}
[c]{|c|}%
A\omega_{1}+B\omega_{2}\\
C\omega_{1}+B\omega_{2}%
\end{array}
\]

\noindent Under such a realization, since this map must restrict to an
automorphism of the subgroup $C_{n_{2}}\times C_{n_{2}}$, the determinant
$\det{(\rho)}=AD-BC$ must be a unit in $C_{n_{2}}$. Assuming $\theta$ has this
form, calculating the two cup products and simplifying, we get
\[
(\theta\circ\eta_{N})\cup(\theta\circ\eta_{N})=N\det{(\theta)}(\omega
_{1}\otimes\omega_{2}-\omega_{2}\otimes\omega_{1}))
\]
and
\[
\eta_{M}\cup\eta_{M}=M(\omega_{1}\otimes\omega_{2}-\omega_{2}\otimes\omega
_{1}).
\]
Since $\det{(\theta)}$ is a unit in $C_{n_{2}}$, if $\gcd{(n,N)}\neq
\gcd{(n,M)}$, then we can never have $N\det{(\theta)}=M$ (else they would be
associates in the ring $C_{n_{2}}$ and thus we would have $\gcd{(n,N)}%
=\gcd{(n,M)}$). Thus if $\gcd{(n,N)}\neq\gcd{(n,M)}$ then $\eta_{N}$ and
$\eta_{M}$ define inequivalent classes. To finish, we need to show if
$\gcd{(n,N)}=\gcd{(n,M)}$, then $\eta_{N}$ and $\eta_{M}$ are equivalent.

Suppose $n=p_{1}\dots p_{s} q_{1} \dots q_{r}$ where the $p_{i}$ and $q_{i}$
are all distinct primes and $d=p_{1} \dots p_{s} =\gcd{(n,N)}$. Observe that
by applying $\mathfrak{S}_{2}$ followed by $\mathfrak{Z}_{1}^{k}$ for any
integer $k$ to $\eta_{N}$, we have the epimorphism $(\alpha_{1} ,\beta_{1}
,\alpha_{2} ,\beta_{2} )\rightarrow(\omega_{1} ,N\omega_{2} ,0,(1-kN)
\omega_{2} )$. Moreover, if $(1-kN,n_{2})=1$, then there exists an integer $a$
which induces an automorphism of $C_{n_{2}}$ via multiplication, such that
$a(1-kN)=1 \mod{(n_{2})}$. Applying this automorphism and reapplying
$\mathfrak{S}_{2}$, we see that $\eta_{N}$ is equivalent to $\eta_{aN}$ for
any such $a$ and $(n_{2},aN)=d$ (since $a$ is a unit in $C_{n_{2}})$.

Assume $((1-kN),n_{2})=1$ and $((1-cN),n_{2})=1$, let $a=(1-kN)^{-1}$ and
$b=(1-cN)^{-1}$ (all taken $\mod{(n_{2})}$). In order to avoid redundancy, we
may assume $0\leq k,c\leq q_{1} \dots q_{r}$. Now if $aN= bN \mod{(n_{2})}$,
then $(1-kN)=(1-cN)$ and consequently $q_{1}\dots q_{r} |(k-c)$. However, this
implies $k=c+aq_{1} \dots q_{r}$ which contradicts that $0\leq k,c\leq q_{1}
\dots q_{r}$ unless $a=0$. Thus $\eta_{aN} =\eta_{bN}$ only when $a=b$.
Therefore, since there are precisely $\phi(n_{2} /d) =(q_{1}-1)(q_{2}%
-1)\dots(q_{r}-1)$ different integers $M$ with $(n_{2},M)=d$, it suffices to
show that there are precisely $(q_{1}-1)(q_{2}-1)\dots(q_{r}-1)$ choices for
$k$ such that $(1-kN,n_{2})=1$.

Observe that $p_{i}$ does not divide $1-kN$ for any value of $i$ or $k$, so it
follows that $(1-kN,n_{2})=1$ provided $T=q_{i_{1}} \dots q_{i_{t}}$ does not
divide $1-kN$ for any divisor $T$ of $q_{1} \dots q_{r}$. Suppose $T=q_{i_{1}}
\dots q_{i_{t}}$ is a divisor of $q_{1}\dots q_{r}$. Then $(T,N)=1$ (since
$\gcd{(n_{2},N)}=d$), so there exists $k$ and $c$ such that $cT+kN=1$, or
$cT=1-kN$. Moreover, if $c_{1}T+k_{1}N=1$ for any other $c_{1}$ and $k_{1}$,
we must have $c_{1} =c+vN$ for some integer $v$ and for every integer $v$, if
we define $c_{1}=c+vN$ and $k_{1} =k-Tv$, then $c_{1}T=1-k_{1}N$. This implies
there will be precisely $q_{1}\dots q_{r}/T$ values of $k$ for which $T$
divides $1-kN$. Since there are $q_{1}\dots q_{r}$ total choices for $k$ and
$q_{1}\dots q_{r}/T$ choices which are divisible by $T$ for each divisor $T$
of $q_{1}\dots q_{r}$, we can form a weighted sum over all the divisors of
$q_{1}\dots q_{r}$ to determine precisely how many are not divisible by any
such $T$. Specifically, the number of values of $k$ with $0\leq k<q_{1} \dots
q_{r}$ such that $(1-kN,q_{1}\dots q_{r})=1$ is
\[
\sum_{n=0}^{r} \bigg[ (-1)^{n} \sum_{r\geq i_{1}> i_{2} > \dots>i_{n}\geq1}
\frac{q_{1}\dots q_{r}}{q_{i_{1}}\dots q_{i_{n}}} \bigg]=(q_{1}-1)(q_{2}%
-1)\dots(q_{r}-1).
\]
The result follows.
\end{proof}

In fact Lemma \ref{lem-equivclasses} can be generalized further so for $p=2$
we don't have $16|n_{2}$ (though the proof is more technical). The next result
is a generalization of the elementary abelian case and the proof is identical.

\begin{lemma}
\label{lem-n2aprime} If $n_{2}=p$ a prime, then there exists $K$ with $r/2\leq
K\leq\mathrm{{min}(}$\textrm{$\rho$}$\mathrm{,r)}$ such that any epimorphism
$\eta\colon\Gamma\rightarrow G$ is equivalent to one and only one of those
described below.%

\[
\eta:=%
\begin{cases}
\eta(\alpha_{i} )=\omega_{i} & \text{$i \leq K$}\\
\eta(\alpha_{i} )=0 & \text{$i > K$}\\
\eta(\beta_{i} )=\omega_{i+K} & \text{$i \leq r-K$}\\
\eta(\beta_{i} )=0 & \text{$i > r-K$}\\
\end{cases}
\]

\end{lemma}

We now have enough information to complete the classification of fixed point
free abelian actions up to genus $65$.

\begin{example}
\label{ex-genus65} First we eliminate all cases which can be classified by our
results. Suppose $G$ has invariant factors $n_{1},n_{2},\dots,n_{r}$. By
Theorem \ref{thm-unramified}, we may assume that $n_{1}\neq p$ for $p$ a prime
and by Lemma \ref{lem-n2aprime}, if the $p$-rank is greater than $1$, we may
assume that $n_{2}\neq p$ for $p$ a prime either. Consequently, if $G$ has
$p$-rank at least $3$ and $n_{3}=a$, then $n_{2}=ab$ and $n_{1}=abc$ where if $a$
is prime, then $b\neq1$. Thus we would have $|G|=a^{3}b^{2}c\geq64$ unless
$a=2$, $b=2$ and $c=1$. However, in this case the invariant factors are
$4,4,2$ and we classified all such epimorphisms in Example \ref{ex-genus33}.
Thus we may assume $G$ has $p$-rank $2$.

For $p$-rank $2$, Lemma \ref{lem-equivclasses} implies we may assume there exists
a prime $p$ such that $p^{2}|n_{2}$ (and consequently $p^{2}|n_{1}$). If
$p\geq3$, then $|G|\geq3^{4}$ and so the genus of the resulting surface will
be at least $82$, so we may assume $p=2$. This means $n_{1}=4k$ for some
integer $k$. If $k>3$, the resulting surface will have genus greater than
$65$, so we only need $n_{1}=n_{2}=4$,  $n_{1}=8$ and $n_{2}=4$, and $n_{1}=12$ and $n_{2}=4$. However, in
all these cases, there is a unique value $N$ for which $(N,n_{2})>1$, so
using the cohomology arguments of Lemma \ref{lem-equivclasses} and the general
form given in Proposition \ref{prop-abclasses}, it is easy to check that there
are just three classes:
\[%
\begin{array}
[c]{l}%
\eta_{1}\colon(\alpha_{1},\beta_{1},\alpha_{2},\beta_{2})\rightarrow
(\omega_{1},\omega_{2},0,0)\\
\eta_{2}\colon(\alpha_{1},\beta_{1},\alpha_{2},\beta_{2})\rightarrow
(\omega_{1},2\omega_{2},\omega_{2},0)\\
\eta_{3}\colon(\alpha_{1},\beta_{1},\alpha_{2},\beta_{2})\rightarrow
(\omega_{1},0,\omega_{2},0)
\end{array}
\]

\end{example}

\section{The totally ramified case\label{Sec-Ramified}\bigskip}

Assume that $G=\mathbb{F}_{p}^{v}$ and that our action is purely ramified,
namely $S/G$ has genus $0,$ and that there are $r$ branch points. We are going
to describe a method for computing the number of equivalence classes of
actions of $G$. In each case there will be several small exceptional primes
and a finite number of infinite families for which the number of actions is
expressed as a polynomial in the prime. The infinite families are defined by
the existence or non-existence of certain roots of unity and so each family is
typically the set of primes in an arithmetic progression. Unfortunately, our
method depends on the enumeration of the finite subgroups of the symmetric
group on $r$ points, so exact general results are impossible for all $p$-ranks and
numbers of branch points. Thus, we will describe the general method, but only
derive the precise details for 3 and 4 branch points (Table 5.1).

We may associate a generating set $\{C_{j}:1\leq j\leq r\}$ with a matrix $X$
\[
\{C_{j}:1\leq j\leq r\}\leftrightarrow\left[
\begin{array}
[c]{cccc}%
X_{1} & X_{2} & \cdots & X_{r}%
\end{array}
\right]  =X
\]
where
\begin{equation}
X\mathrm{\ is\ a\ }v\times r\mathrm{\ matrix\ of\ rank\ }v \label{eq-Xcond}%
\end{equation}
such that%
\begin{equation}
X_{1}+X_{2}+\cdots+X_{r}=0. \label{eq-Xsumzero}%
\end{equation}
Note that $1\leq v<r$ and $r\geq2.$ Let us explicitly define the action as we
shall need it in our calculations. Let $g\in GL(v,\mathbb{F}_{p})$ be a
$v\times v$ invertible matrix over $\mathbb{F}_{p},$ and $\alpha\in\Sigma
_{r}.$ Let $\pi_{\alpha}$ be the standard permutation matrix defined by
\[
\pi_{\alpha}=\left[
\begin{array}
[c]{c}%
E_{\alpha(1)}\\
E_{\alpha(2)}\\
\vdots\\
E_{\alpha(r)}%
\end{array}
\right]
\]
where $E_{1},\ldots,E_{n}$ are the rows of the identity matrix$.$ Then
$\pi_{\alpha\beta}=$ $\pi_{\alpha}\pi_{\beta}$ and the action $\mathrm{Aut}%
(G)\times\mathrm{Aut}(\Gamma)$ on the set of matrices $X$ is the action of
$GL(v,\mathbb{F}_{p})\times\Sigma_{r}$ given by
\begin{equation}
(g,\alpha)\cdot X=gX\pi_{\alpha}^{\top}. \label{eq-MCGaction}%
\end{equation}
i.e.,%
\[
gX\pi_{\alpha}^{\top}=\left[
\begin{array}
[c]{cccc}%
gX_{\alpha\left(  1\right)  } & gX_{\alpha\left(  2\right)  } & \cdots &
gX_{\alpha\left(  r\right)  }%
\end{array}
\right]
\]
The problem of computing the number of equivalence classes is the same as
computing the number of orbits under the given action. Before we start setting
up the machinery for the computation of the orbits, let us prove a trivial
case for all $p$-ranks.

\begin{proposition}
\label{pr-onemore}Suppose that the number of branch points is one more than
the $p$-rank of $G.$ Then there is a unique action of $G$.
\end{proposition}

\begin{proof}
Let $g=\left[
\begin{array}
[c]{cccc}%
X_{1} & X_{2} & \cdots & X_{v}%
\end{array}
\right]  $ be the square matrix obtained by dropping the last column of
$X$. Since $X_{1}+X_{2}+\cdots+X_{r}=0$ then the column space of $g$ is the
same as the column space of $X.$ It follows the that $g$ is full rank and
hence invertible. Then $g^{-1}X=\left[
\begin{array}
[c]{cc}%
I_{v} & -E_{v}%
\end{array}
\right]  $ where $I_{r}$ is the $v\times v$ identity matrix and $E_{v}$ a
column of ones.
\end{proof}

\subsection{Orbit spaces via equisymmetry\label{subsec-eqsym}}

Let $\Omega$ be a finite set upon which the finite group $G$ (different from
our $G$ above) acts. Consider the following standard constructs. Given a
subset $Y\subseteq\Omega$ or $H\subseteq G$ let
\[
G_{Y}=\left\{  g\in G:gx=x,\forall x\in Y\right\}  \text{ and }G_{x}%
=G_{\left\{  x\right\}  }%
\]
and
\[
\Omega^{H}=\left\{  x\in\Omega:gx=x,\forall g\in H\right\}  \text{ and }%
\Omega^{g}=\Omega^{\left\{  g\right\}  }%
\]
For $H\subseteq G$ let $\overline{H}$ denote the \emph{action closure} of $H,$
the largest group fixing all the points fixed by all elements of $H,$
\[
\overline{H}=G_{\Omega^{H}}=\bigcap\limits_{x\in\Omega^{H}}G_{x}.
\]
Also $\overline{Y}$ denotes the \emph{action closure}\ of a subset $Y$
$\subseteq\Omega,$ the set of all points fixed by all elements of $G$ fixing
all of $Y$,
\[
\overline{Y}=\Omega^{G_{Y}}=\bigcap\limits_{g\in G_{Y}}\Omega^{g}.
\]
Next for a subgroup $H\subseteq G$ let $\left\langle H\right\rangle $ denote
the conjugacy class of subgroups of $G$ determined by $H.$ The set of
conjugacy classes has a partial order $\left\langle H_{1}\right\rangle
\prec\left\langle H_{2}\right\rangle $ if and only if $H_{1}\subset
gH_{2}g^{-1}$ for some $g\in G.$ We make the following additional definitions.

\begin{definition}
Two points $x,y\in\Omega$ are called \emph{equi-isotropic }if\emph{\ }%
$G_{x}=G_{y}$ and \emph{equisymmetric} if $G_{x}$ and $G_{y}$ are conjugate
subgroups. The group $G_{x}$ is called the \emph{isotropy type} of $x$ and
$\left\langle G_{x}\right\rangle $ is called the \emph{symmetry type} of $x$.
If $H=G_{x}$ for some $x$ then $\mathring{\Omega}^{H},$ the $H$%
\emph{-isotropic stratum} or the \emph{isotropic stratum of }$x$\emph{\ }is
the set of all points with isotropy type $H.$
\[
\mathring{\Omega}^{H}=\left\{  y\in\Omega:G_{y}=H\right\}
\]
Likewise if $H=G_{x}$ for some $x$ then $\left\langle \mathring{\Omega}%
^{H}\right\rangle ,$ the $H$\emph{-equisymmetric\ stratum} or the
\emph{equisymmetric\ stratum of }$x,$\emph{\ }is the set of all points with
symmetry type $\left\langle H\right\rangle .$
\[
\left\langle \mathring{\Omega}^{H}\right\rangle =\left\{  y\in\Omega:G_{y}%
\in\left\langle H\right\rangle \right\}
\]

\end{definition}

The action closure\ $\Omega^{H}$ of $\mathring{\Omega}^{H}$ is the set of
points with $\ H$-isotropy or greater, i.e., $\Omega^{H}=\left\{  x\in
\Omega:H\subseteq G_{x}\right\}  .$ It is easily seen that
\begin{equation}
\Omega^{H}=\mathring{\Omega}^{H}\cup\bigcup\limits_{K}\mathring{\Omega}^{K}
\label{eq-isotropy1}%
\end{equation}
a disjoint union for a set of closed groups $K$ strictly containing $K\supset
H.\ $We can actually take the union over all subgroups $K\supset H$ fixing any
point, though some of the $\mathring{\Omega}^{K}$ may be empty. We get a
similar union
\begin{equation}
\left\langle \Omega^{H}\right\rangle =\left\langle \mathring{\Omega}%
^{H}\right\rangle \cup\bigcup\limits_{K}\left\langle \mathring{\Omega}%
^{K}\right\rangle , \label{eq-symmetry1}%
\end{equation}
for a set of closed point stabilizers $K$ such that $\left\langle
K\right\rangle \succ\left\langle H\right\rangle $. $\ $ The equations
(\ref{eq-isotropy1}) and (\ref{eq-symmetry1}) may be rewritten to iteratively
compute
\begin{equation}
\left\vert \mathring{\Omega}^{H}\right\vert =\left\vert \Omega^{H}\right\vert
-\sum\limits_{K}\left\vert \mathring{\Omega}^{K}\right\vert
\label{eq-isotropy2}%
\end{equation}%
\begin{equation}
\left\vert \left\langle \mathring{\Omega}^{H}\right\rangle \right\vert
=\left\vert \left\langle \Omega^{H}\right\rangle \right\vert -\sum
\limits_{K}\left\vert \left\langle \mathring{\Omega}^{K}\right\rangle
\right\vert \label{eq-symmetry2}%
\end{equation}
where the sums are over the same set of subgroups as above.

If $x$ and $y$ have the same symmetry type, i.e., $G_{x}$ and $G_{y}$ are
conjugate, then the orbits $Gx$ and $Gy$ have the same size $\left\vert
G\right\vert /\left\vert G_{x}\right\vert $. Since $G_{gx}=gG_{x}g^{-1},$ then
$g\mathring{\Omega}^{H}=\mathring{\Omega}^{gHg^{-1}}$ for a point stabilizer
$H$ and so the set $\left\langle \mathring{\Omega}^{H}\right\rangle $ is
$G$-invariant. Indeed we have a disjoint union%
\[
\left\langle \mathring{\Omega}^{H}\right\rangle =\bigcup\limits_{K\in
\left\langle H\right\rangle }\mathring{\Omega}^{K}=\bigcup\limits_{g\in
G/N_{G}(H)}g\mathring{\Omega}^{H}%
\]
from which we get $\left\vert \left\langle \mathring{\Omega}^{H}\right\rangle
\right\vert =$ $\left\vert \left\langle H\right\rangle \right\vert \left\vert
\mathring{\Omega}^{H}\right\vert .$ Since $\left\langle \mathring{\Omega}%
^{H}\right\rangle \rightarrow\left\langle \mathring{\Omega}^{H}\right\rangle
/G$ is a $\left\vert G\right\vert /\left\vert H\right\vert $ to $1$ map
\begin{equation}
\left\vert \left\langle \mathring{\Omega}^{H}\right\rangle /G\right\vert
=\frac{\left\vert H\right\vert }{\left\vert G\right\vert }\left\vert
\left\langle \mathring{\Omega}^{H}\right\rangle \right\vert =\frac{\left\vert
H\right\vert }{\left\vert N_{G}(H)\right\vert }\left\vert \mathring{\Omega
}^{H}\right\vert \label{eq-eqsymandeqiso}%
\end{equation}
and it follows that
\begin{equation}
\left\vert \Omega/G\right\vert =\sum\limits_{H}\left\vert \left\langle
\mathring{\Omega}^{H}\right\rangle /G\right\vert =\sum\limits_{H}%
\frac{\left\vert H\right\vert }{\left\vert N_{G}(H)\right\vert }\left\vert
\mathring{\Omega}^{H}\right\vert , \label{eq-eqsymsum}%
\end{equation}
where $\left\langle H\right\rangle $ runs over all conjugacy classes of point
stabilizers. So computing $\left\vert \Omega/G\right\vert $ is just a
matter of determining the point stabilizers $H$ and the equisymmetric or the
equi-isotropic strata. The cardinalities of the closed subsets $\Omega^{H}$
are generally easier to calculate directly which is why we use the above formulas.

\subsection{Matrix formulation of M\"{o}bius inversion \label{subsec-matrix}}

The iterative calculation of $\left\vert \Omega/G\right\vert $ via the
formulas (\ref{eq-isotropy2}) and (\ref{eq-eqsymsum}), which amounts to
M\"{o}bius inversion, can be put into a very succinct matrix format. Consider
the set of subgroups $H\subseteq G$ such that $\Omega^{H}$ is non-empty. For
each conjugacy class of subgroups in this set select a representative and then
list the subgroups $\left\langle id\right\rangle =H_{1},H_{2},\ldots,H_{s},$
in such a way that $\left\langle H_{i}\right\rangle \subseteq$ $\left\langle
H_{j}\right\rangle $ $\Longrightarrow i\leq j$. This can always be done, though
generally not in a unique way. From equation (\ref{eq-isotropy1}) we have
\begin{align}
\left\vert \Omega^{H_{i}}\right\vert  &  =\left\vert \mathring{\Omega}^{H_{i}%
}\right\vert +\sum_{K\supset H_{i}}\left\vert \mathring{\Omega}^{K}\right\vert
\label{eq-matmob1}\\
&  =\left\vert \mathring{\Omega}^{H_{i}}\right\vert +\sum_{j>i}\left(
\sum_{K\in\left\langle H_{j}\right\rangle ,H_{i}\subset K}\left\vert
\mathring{\Omega}^{K}\right\vert \right) \nonumber
\end{align}
For $H,K\subseteq G,$ let%
\begin{align}
U(H,K)  &  =\{gKg^{-1}:g\in G,H\subseteq gKg^{-1}\}\label{eq-UHK}\\
D(H,K)  &  =\{g^{-1}Hg:g\in G,g^{-1}Hg\subseteq K\} \label{eq-DHK}%
\end{align}
and set $u_{i,j}=\left\vert U(H_{i},H_{j})\right\vert ,\ d_{i,j}=\left\vert
D(H_{i},H_{j})\right\vert .$ Observe that $u_{i,i}=d_{i,i}=1$ and
$u_{i,j}=d_{i,j}=0$ if $i>j$ and from Lemma \ref{lm-udcount} below we have
\begin{equation}
u_{i,j}=d_{i,j}\frac{\left\vert \left\langle H_{j}\right\rangle \right\vert
}{\left\vert \left\langle H_{i}\right\rangle \right\vert }.
\label{eq-updowncount1}%
\end{equation}
Then equation (\ref{eq-matmob1}) can be rewritten
\begin{align}
\left\vert \Omega^{H_{i}}\right\vert  &  =\left\vert \mathring{\Omega}^{H_{i}%
}\right\vert +\sum_{i<j}u_{i,j}\left\vert \mathring{\Omega}^{H_{j}}\right\vert
\label{eq-matmob3}\\
&  =\left\vert \mathring{\Omega}^{H_{i}}\right\vert +\sum_{i<j}d_{i,j}%
\frac{\left\vert \left\langle H_{j}\right\rangle \right\vert }{\left\vert
\left\langle H_{i}\right\rangle \right\vert }\left\vert \mathring{\Omega
}^{H_{j}}\right\vert \nonumber
\end{align}
and
\begin{equation}
\left\vert \left\langle \mathring{\Omega}^{H_{i}}\right\rangle \right\vert
=\left\vert \mathring{\Omega}^{H_{i}}\right\vert \left\vert \left\langle
H_{i}\right\rangle \right\vert . \label{eq-matmob5}%
\end{equation}
Now let us convert to matrix notation. Define
\begin{align*}
U  &  =\left[  u_{i,j}\right]  ,\text{ }D=\left[  d_{i,j}\right] \\
S  &  =\mathrm{diag}(\left\vert \left\langle H_{1}\right\rangle \right\vert
,\left\vert \left\langle H_{2}\right\rangle \right\vert ,\ldots,\left\vert
\left\langle H_{s}\right\rangle \right\vert )\\
N  &  =\mathrm{diag}(\left\vert N_{G}\left(  H_{1}\right)  \right\vert
,\left\vert N_{G}\left(  H_{2}\right)  \right\vert ,\ldots,\left\vert
N_{G}\left(  H_{2}\right)  \right\vert )\\
T  &  =\mathrm{diag}(\left\vert H_{1}\right\vert ,\left\vert H_{2}\right\vert
,\ldots,\left\vert H_{s}\right\vert ).
\end{align*}
The matrices $U=\left[  u_{i,j}\right]  $ and $D=\left[  d_{i,j}\right]  $ are
upper triangular with $1$'s on the diagonal. It follows from equation
\ref{eq-updowncount1} that $U=S^{-1}DS$. Let $L,$ $L^{\circ},E^{\circ
},O^{\circ}$ be the vectors
\begin{align*}
L  &  =\left[
\begin{array}
[c]{cccc}%
\left\vert \Omega^{H_{1}}\right\vert  & \left\vert \Omega^{H_{2}}\right\vert
& \cdots & \left\vert \Omega^{H_{s}}\right\vert
\end{array}
\right]  ^{\top}\\
L^{\circ}  &  =\left[
\begin{array}
[c]{cccc}%
\left\vert \mathring{\Omega}^{H_{1}}\right\vert  & \left\vert \mathring
{\Omega}^{H_{2}}\right\vert  & \cdots & \left\vert \mathring{\Omega}^{H_{s}%
}\right\vert
\end{array}
\right]  ^{\top}\\
E^{\circ}  &  =\left[
\begin{array}
[c]{cccc}%
\left\vert \left\langle \mathring{\Omega}^{H_{1}}\right\rangle \right\vert  &
\left\vert \left\langle \mathring{\Omega}^{H_{2}}\right\rangle \right\vert  &
\cdots & \left\vert \left\langle \mathring{\Omega}^{H_{s}}\right\rangle
\right\vert
\end{array}
\right]  ^{\top}\\
O^{\circ}  &  =\left[
\begin{array}
[c]{cccc}%
\left\vert \left\langle \mathring{\Omega}^{H_{1}}\right\rangle /G\right\vert
& \left\vert \left\langle \mathring{\Omega}^{H_{2}}\right\rangle /G\right\vert
& \cdots & \left\vert \left\langle \mathring{\Omega}^{H_{s}}\right\rangle
/G\right\vert
\end{array}
\right]  ^{\top}.
\end{align*}
The vectors $E^{\circ}$ and $O^{\circ}$ are the cardinalities of the various
equisymmetric strata and corresponding orbits, $\left\vert \Omega/G\right\vert
$ is simply the sum of the entries of $O^{\circ}.$ The equations
(\ref{eq-matmob1}) to (\ref{eq-matmob5}) can be translated into the following
matrix equations:
\begin{align*}
L  &  =UL^{\circ}=S^{-1}DSL^{\circ}\\
L^{\circ}  &  =U^{-1}L=S^{-1}D^{-1}SL\\
E^{\circ}  &  =SL^{\circ}=SS^{-1}D^{-1}SL=D^{-1}SL\\
O^{\circ}  &  =\frac{1}{\left\vert G\right\vert }TE^{\circ}=\frac
{1}{\left\vert G\right\vert }TD^{-1}SL.
\end{align*}
Finally we have the matrix equation $N^{-1}=\frac{1}{\left\vert G\right\vert
}S$ relating normalizers and conjugacy classes resulting in our final
equation
\begin{equation}
O^{\circ}=TD^{-1}N^{-1}L. \label{eq-orbcount}%
\end{equation}

The following lemma relates the $u_{i,j}$ and $d_{i,j}.$

\begin{lemma}
\label{lm-udcount}Let $H\subseteq K$ be subgroups of a finite group $G$. Let
$U(H,K)$ and $D(H,K)$ be the sets defined in equations (\ref{eq-UHK}) and
(\ref{eq-DHK}). Then%
\begin{equation}
U(H,K)\left\vert \left\langle H\right\rangle \right\vert =\left\vert
D(H,K)\right\vert \left\vert \left\langle K\right\rangle \right\vert .
\label{eq-updowncount2}%
\end{equation}

\end{lemma}

\begin{proof}
Let $P=\left\{  (S,T)\in\left\langle H\right\rangle \times\left\langle
K\right\rangle :S\subseteq T\right\}  .$ By considering the projections
$P\rightarrow\left\langle H\right\rangle ,$ $P\rightarrow\left\langle
K\right\rangle $ we see that $\left\vert U(H,K)\right\vert \left\vert
\left\langle H\right\rangle \right\vert $ $=$ $\left\vert P\right\vert $ $=$
$\left\vert D(H,K)\right\vert \left\vert \left\langle K\right\rangle
\right\vert .$
\end{proof}

\begin{remark}
The entries of $D$ can be found by enumerating the subgroups of $H_{j}$ and
then determining to which $H_{i}$ they are conjugate. In our calculation we
will just work with subgroups of small symmetric groups.
\end{remark}

\subsection{Structure of stabilizers\label{subsec-stabstruct}}

Let $\Omega$ be the set of matrices $X$ satisfying conditions (\ref{eq-Xcond})
and (\ref{eq-Xsumzero}), and $G=GL(v,\mathbb{F}_{p})\times\Sigma_{r}$ with the
action given in equation (\ref{eq-MCGaction}). There are two canonical
representations $\rho_{1}:G\rightarrow GL(v,\mathbb{F}_{p}),$ $(g,\alpha
)\rightarrow g$ and $\rho_{2}:G\rightarrow\Sigma_{r}\rightarrow
GL(r,\mathbb{F}_{p}),$ $(g,\alpha)\rightarrow\pi_{\alpha}.$ We are going to
describe the structure of stabilizers $H\subseteq G$ in terms of the two
representations. Let $X_{0}\in\Omega$ and suppose that $H$ is a subgroup
fixing $X_{0},$ (the subgroup need not be action closed). If $X$ is any matrix
stabilized by $H$ then $\rho_{1}(h)X\rho_{2}^{\top}(h)=X,$ $h\in H$ or
\begin{equation}
\rho_{1}(h)X=X\rho_{2}(h),h\in H, \label{eq-intertwine1}%
\end{equation}
i.e., $X\in$ \textrm{Hom}$_{H}(\mathbb{F}_{p}^{r},\mathbb{F}_{p}^{v})$ is an
intertwining operator. Now, the restricted map $p_{2}:H\rightarrow\Sigma_{r},$
is 1-1, for if $(g,\alpha)\in$ $\mathrm{ker}(p_{2})$ then $\alpha=id$ and
$X_{0}=gX_{0}.$ Since $X_{0}$ has rank $v$ then $g=I.$ Let $H^{\prime}%
=p_{2}(H)$ denote the isomorphic image of $H$ in $\Sigma_{r}.$ Reversing the
process, we can construct candidate stabilizers by selecting a subgroup
$H^{\prime}\subseteq\Sigma_{r},$ and a suitable representation $q=\rho
_{1}\circ p_{2}^{-1}:H^{\prime}\rightarrow GL(v,\mathbb{F}_{p})$ and then
construct $H$ by
\begin{equation}
H=\{(q(\alpha),\alpha):\alpha\in H^{\prime}\} \label{eq-stabconstruct}%
\end{equation}
It is not clear that an arbitrary $H$ so constructed fixes any matrices in
$\Omega,$ in fact the representation $q$ must satisfy certain properties which
we proceed to discuss.\bigskip

For the remainder of our discussion we are going to assume that $p\nmid
\left\vert H\right\vert ,$ so that we may use the theory of reducible
representations to aid our computations. Since $\left\vert H\right\vert $
divides $r!$ we are only excluding the small primes $p\leq r.$ Both
$H$-modules $\mathbb{F}_{p}^{v},$ $\mathbb{F}_{p}^{r}$ can be rewritten as a
direct sum of $\mathbb{F}_{p}$-irreducibles, $\mathbb{F}_{p}^{v}=V_{1}\oplus
V_{2}\oplus\cdots\oplus V_{m},$ $\mathbb{F}_{p}^{r}=W_{1}\oplus W_{2}%
\oplus\cdots\oplus W_{n}.$ By a change in coordinates, we may write
\begin{align}
A^{-1}\rho_{1}(h)A  &  =\left[
\begin{array}
[c]{cccc}%
\phi_{1}(h) & 0 & \cdots & 0\\
0 & \phi_{2}(h) & \cdots & 0\\
\vdots & \vdots & \ddots & \vdots\\
0 & 0 & \cdots & \phi_{m}(h)
\end{array}
\right] \label{eq-diag}\\
B^{-1}\rho_{2}(h)B  &  =\left[
\begin{array}
[c]{cccc}%
\psi_{1}(h) & 0 & \cdots & 0\\
0 & \psi_{2}(h) & \cdots & 0\\
\vdots & \vdots & \ddots & \vdots\\
0 & 0 & \cdots & \psi_{n}(h)
\end{array}
\right] \nonumber
\end{align}
where each $\phi_{i},$ and $\psi_{j}$ is irreducible on $V_{i}$ and $W_{j}$
respectively. Writing $Y=A^{-1}XB$ in corresponding block form%
\begin{equation}
Y=A^{-1}XB=\left[
\begin{array}
[c]{cccc}%
Y_{1,1} & Y_{1,2} & \cdots & Y_{1,n}\\
Y_{2,2} & Y_{2,2} & \cdots & Y_{2,n}\\
\vdots & \vdots & \ddots & \vdots\\
Y_{m,1} & Y_{m,2} & \cdots & Y_{m,n}%
\end{array}
\right]  \label{eq-Ydecomp}%
\end{equation}
the intertwining property has an equivalent formulation
\begin{equation}
\rho_{1}(h)X=X\rho_{2}(h)\Longleftrightarrow\forall_{i,j}\text{ }\phi
_{i}(h)Y_{i,j}=Y_{i,j}\psi_{j}(h), \label{eq-intertwine2}%
\end{equation}
and hence we may consider $Y_{i,j}\in\mathrm{Hom}_{H}(W_{j},V_{i}).$ Now, as
$H$ acts irreducibly on $V_{i}$ and $W_{j},$ then $Y_{i,j}=0$ if $\phi_{i},$
$\psi_{j}$ are inequivalent. If $\phi_{i}$ and $\psi_{j}$ are equivalent then
$Y_{i,j}$ is invertible, if non-zero, by Schur's Lemma. The number of
invertible $Y_{i.j}$ is $p^{k}-1,$ where $k$ is the number of irreducibles
into which $\phi_{i}=$ $\psi_{j}$ splits in the algebraic closure of
$\mathbb{F}_{p}.$ We choose $\psi_{n}$ to be the trivial representation acting
on the space $W_{n}=\left\langle \left[
\begin{array}
[c]{cccc}%
1 & 1 & \ldots & 1
\end{array}
\right]  ^{\top}\right\rangle .$ Since $X\left[
\begin{array}
[c]{cccc}%
1 & 1 & \ldots & 1
\end{array}
\right]  ^{\top}=0$ then $Y_{i,n}=0$ for all $i$. For convenience we introduce
a special name for the representation on $W_{1}\oplus\cdots\oplus W_{n-1}.$

\begin{definition}
For a subgroup $H^{\prime}\subseteq$ $\Sigma_{r}$ acting naturally on the
space $\mathbb{F}_{p}^{r}$ write $\mathbb{F}_{p}^{r}=W\oplus$ $\left\langle
\left[
\begin{array}
[c]{cccc}%
1 & 1 & \ldots & 1
\end{array}
\right]  ^{\top}\right\rangle $ as a direct sum of subrepresentations. If
$\pi$ denotes the natural representation of $H^{\prime},$ we call $\pi_{W},$
the natural representation restricted to $W,$ the reduced natural
representation of $H^{\prime}$ and denote it by $\overline{\pi}.$ If $p\nmid
r$ then may take $W=$ $\left\{  (x_{1},\ldots,x_{r})\in\mathbb{F}_{p}%
^{r}:x_{1}+\cdots+x_{r}=0\right\}  .$
\end{definition}

\begin{remark}
Though a bit pedantic, we remind ourselves on how to construct $A$ and $B$
since the details will be used in a later proof. The matrix $A$ may be
constructed by choosing column vectors that span $V_{1}$ then adding vectors
that span $V_{2}$ and so on to get $A=\left[
\begin{array}
[c]{cccc}%
A_{1} & A_{2} & \cdots & A_{m}%
\end{array}
\right]  $, where $\ A_{i}$ is a $v\times v_{i}$ matrix. Construct $B=\left[
\begin{array}
[c]{cccc}%
B_{1} & B_{2} & \cdots & B_{n}%
\end{array}
\right]  ,$ similarly. Now write $A^{-1}=\left[
\begin{array}
[c]{cccc}%
\left(  A_{1}^{\ast}\right)  ^{\top} & \left(  A_{2}^{\ast}\right)  ^{\top} &
\cdots & \left(  A_{m}^{\ast}\right)  ^{\top}%
\end{array}
\right]  ^{\top}$ as a column of matrices with $A_{i}^{\ast}$ of size
$v_{i}\times v.$ Define the $B_{j}^{\ast}$ similarly. We have the following
transformation formula
\begin{equation}
Y_{i,j}=A_{i}^{\ast}XB_{j},\text{ }\mathrm{\ }X=\sum\limits_{i,j}A_{i}%
Y_{i,j}B_{j}^{\ast} \label{eq-intertwine3}%
\end{equation}

\end{remark}

The next definition helps us characterize the representations we are looking
for in the defining equation (\ref{eq-stabconstruct}).

\begin{definition}
\label{def-dom}Let $H$ be a finite group and $k$ a field whose characteristic
is coprime to $\left\vert H\right\vert $. Let $\rho_{1},$ $\rho_{2}$ be two
representations of finite degree on $k$ vector spaces $V,W$. We say that
$\rho_{2}$ dominates $\rho_{1}$ $(\rho_{1}\preceq\rho_{2}) $ if and only if
for each $k$-irreducible representation $\phi$ of $H, $ $\left\langle
\phi,\rho_{1}\right\rangle \leq$ $\left\langle \phi,\rho_{2}\right\rangle $
where $\left\langle \phi,\rho\right\rangle $ is the number of times the
representation $\phi$ is contained in a representation $\rho$ of $H.$
\end{definition}

\begin{proposition}
\label{pr-dom}Let $H,$ $k$ and $\rho_{1},$ $\rho_{2}$ be as in Definition
\ref{def-dom}. Then $\rho_{1}\preceq\rho_{2}$ if and only if the space of
intertwining operators $\mathrm{Hom}_{H}(W,V)$ contains a surjective map.
Similarly $\rho_{2}\preceq\rho_{1}$ if and only if $\mathrm{Hom}_{H}(W,V)$
contains an injective map.
\end{proposition}

When restated in terms $H$-modules over $\mathbb{F}_{p}$ the proposition is
self-evident. In proof of Proposition \ref{pr-stabchar}, which characterizes
the representations defining subgroups that fix a point, we will see that the
dominance condition $\rho_{1}\preceq\rho_{2}$ is equivalent to the existence
of an intertwining operator $X$ with linearly independent rows.

\begin{remark}
\label{rk-orbits}Suppose the action of $H$ has several orbits on
$\{1,2,\ldots,r\}$ and that by selecting a suitable conjugate of $H^{\prime}$
each orbit is an interval of $t$ integers $\left\{  s,s+1\ldots,s+t-1\right\}
$. We may write
\begin{equation}
X=\left[
\begin{array}
[c]{cccc}%
Z_{1} & Z_{2} & \cdots & Z_{l}%
\end{array}
\right]  \label{eq-suborbits1}%
\end{equation}
as a block matrix where each block is defined by an orbit and hence invariant
under $H^{\prime}.$ Then $\rho_{1}$ must satisfy
\begin{equation}
\rho_{1}(h)Z_{k}=Z_{k}\rho_{2,k}(h),k=1,\ldots,l \label{eq-suborbits2}%
\end{equation}
where $\rho_{2,k}(h)$ is the induced permutation matrix on the orbit defining
$Z_{k}$. We use these representations to characterize when a pair $H^{\prime
}\subseteq\Sigma_{r},$ $q:H^{\prime}\rightarrow GL(v,\mathbb{F}_{p})$ fixes a
point in $\Omega.$
\end{remark}

\begin{proposition}
\label{pr-stabchar} Assume that $v>1$ if $p=2$. Let $H^{\prime}\subseteq
\Sigma_{r}$ such that $p\nmid\left\vert H\right\vert ,$ let $q:H^{\prime
}\rightarrow GL(v,\mathbb{F}_{p})$ be an arbitrary representation, and let
$H=\{(q(\alpha),\alpha):\alpha\in H^{\prime}\}.$ Let $\rho,\overline{\rho}$ be
the standard and reduced representations of $H^{\prime}$ afforded by
$H^{\prime}\subseteq\Sigma_{r},$ and $\rho_{i}=\rho_{2,i}\circ p_{2}^{-1}$ be
the representations of $H^{\prime}$ determined by the orbits of $H^{\prime}$
on $\{1,2,\ldots,r\}.$ Then $H$ fixes a point in $\Omega$ if and only if all
the following conditions hold.

\begin{itemize}
\item $q$ is dominated by $\overline{\rho}$ $,$ and

\item $q$ and $\rho_{k}$ have a common irreducible for $k=1,\ldots,l$.
\end{itemize}
\end{proposition}

\begin{proof}
By taking a conjugate of $H^{\prime}$ we may assume that the orbits occur in
intervals as Remark \ref{rk-orbits}. Let $\rho_{1}$ and $\rho_{2},$
$\rho_{2,k},$ $Z_{k}$ be as defined previously so that $\rho_{1}=q\circ
p_{2},$ and, $\rho_{2}=\rho\circ p_{2},$ $\overline{\rho_{2}}=\overline{\rho
}\circ p_{2},$ $\rho_{2,i}=\rho_{i}\circ p_{2}.$ To show that $H$ fixes an
element of $\Omega$ we need to find a $v\times r$ matrix $X$ such that

\begin{itemize}
\item $\rho_{1}(h)X=X\rho_{2}(h),h\in H,$

\item $X_{1}+X_{2}+\cdots+X_{r}=0,$

\item $X$ is surjective, i.e., has linearly independent rows,

\item each column $X_{i}$ is non-zero.
\end{itemize}

We leave it to the reader to show that the conditions are necessary, we show
how to construct an $X$ if the conditions hold. Our proof depends on carefully
setting up the diagonalizing matrices $A$ and $B,$ and then carefully
selecting the $Y_{i,j}$. Because $(g,\alpha)$ in $H$ acts on a column $X_{j}$
by $X_{j}\rightarrow gX_{\alpha(j)},$ then a $Z_{k}$ defined by an $H^{\prime
}$ orbit will be non-zero if and only if all the columns of the given $Z_{k}$
are non-zero. Thus we merely need to construct $X$ satisfying the first
three bullets and all $Z_{i}$ non-zero. The diagonalizing matrix $A$ may be
constructed according to any decomposition, we need more care with $B.$
Decompose $\rho_{2,1}$ into $H$-irreducibles and then place the corresponding
basis vectors into $B$ as columns as described in Remark \ref{rk-orbits}. Next
decompose $\rho_{2,2}$ and add the basis vectors as columns, and continue on
to the last orbit. Each orbit determines a unique trivial subrepresentation of
the $\rho_{2,k}.$ A spanning vector for this subrepresentation is the vector
with 1's in the locations corresponding to the orbit and zeros elsewhere. By
construction, for each of these vectors, $B$ contains a column which is a
scalar multiple of this vector. We assume that the scalar is $1$ and that the
last column of $B$ is the vector corresponding to the orbit defining
$\rho_{2,l}.$ Now the sum of these columns is $\left[
\begin{array}
[c]{cccc}%
1 & 1 & \ldots & 1
\end{array}
\right]  ^{\top}.$ If we replace the last column of $B$ with the vector
$\left[
\begin{array}
[c]{cccc}%
1 & 1 & \ldots & 1
\end{array}
\right]  ^{\top}$ we obtain a matrix $C$ with the same column span as $B$ and
hence $C$ is invertible. Moreover as we exchanged one invariant vector for
another $B^{-1}\rho_{2}(h)B=$ $C^{-1}\rho_{2}(h)C.$ Thus we can assume that
$B$ has $B_{n}=\left[
\begin{array}
[c]{cccc}%
1 & 1 & \ldots & 1
\end{array}
\right]  ^{\top}$ as its last column.

Now we are going to further modify $A$ and $B$ as follows. We assume that we
have selected our representations $\phi_{i}$ and $\psi_{j}$ so that if any two are 
equivalent then they are equal. This can be achieved by modifying the columns
of $A$ or $B$ corresponding to the $H$-invariant subspace corresponding to a
$\phi_{i}$ or $\psi_{j}$. Note that we do not need to alter the last column of
$B.$ In this circumstance we define the matrix components of $Y$ in equation
(\ref{eq-Ydecomp})
\begin{align}
Y_{i,j}  &  =y_{i,j}I_{v_{i}}\text{ \textrm{if} }\phi_{i}=\psi_{j}%
\label{eq-intertwine4}\\
Y_{i,j}  &  =0,\text{ }\mathrm{otherwise.}\nonumber
\end{align}
By equation (\ref{eq-intertwine2}) the above equations always define an
intertwining operator constructed via equation (\ref{eq-intertwine3}). Now
$X_{1}+X_{2}+\cdots+X_{r}=XB_{n},$ where $B_{n}=\left[
\begin{array}
[c]{cccc}%
1 & 1 & \ldots & 1
\end{array}
\right]  ^{\top}.$ But $XB_{n}$ equals the last column of $XB=AY.$ It follows
then that $X_{1}+X_{2}+\cdots+X_{r}=0$ if and only if the last column of $Y$
is zero, i.e., $Y_{i.n}=0$ for all $i$. Next we modify $Y$ so that $X$ is
surjective. From the hypotheses $\overline{\rho_{2}}$ dominates $\rho_{1},$
and hence for each irreducible $\phi_{i}$ there can be chosen a $\psi_{j(i)}$
with $\phi_{i}=\psi_{j(i)},$ $j(i)<n,$ and such that for distinct $i_{1}%
,i_{2}$ we have $j(i_{i})\neq$ $j(i_{i}).$ If we select $y_{i,j(i)}\neq0$ and
all other $y_{i,j}=0$ then the resulting matrix $Y,$ and consequently $X$, has
linearly independent rows. (This argument shows surjectivity in Proposition
\ref{pr-dom}.)

Finally, we show that the $Y_{i,j}$ can be chosen so that each $Z_{k}$ in
equation (\ref{eq-suborbits1}) is non-zero. By equation (\ref{eq-suborbits1})
$Y=A^{-1}XB=\left[
\begin{array}
[c]{cccc}%
A^{-1}Z_{1}B & A^{-1}Z_{2}B & \cdots & A^{-1}Z_{l}B
\end{array}
\right]  $ so if some $Z_{i}=0$ then $Y$ has a block column equal to zero and
by construction there is a set $J_{0}$ of $j$'s such that $Y_{i,j}=0,$ $1\leq
i\leq m,$ $j\in J_{0}.$ The $J_{0}$ consists of all $j$ such that $\psi_{j}$
comes from the subspace determined by the orbit corresponding to $Z_{k}.$ By
assumption $\rho_{1}$ and $\rho_{2,k}$ have a common irreducible. Then we can
set some $Y_{i,j}=y_{i,j}I_{v_{i}}$ for $i$ and a $j\in J_{0},$ and any
$y_{i,j}\neq0.$ The modified $Y$ will still define an intertwining operator,
and $Y$ will remain full rank since we are modifying columns that were
initially zero. Thus we can guarantee that all $Z_{k}$ are non-zero, and still
leave $X$ a surjective intertwining operator. The potential problem is if we
were forced to change something in the last column of $Y.$ This can only
happen if the trivial representation is the only common irreducible between
$\rho_{2,l}$ and $\overline{\rho_{2}}.$ Now suppose that there is another
$\rho_{2,l^{\prime}}$ such that $\overline{\rho_{2}}$ and $\rho_{2,l^{\prime}%
}$ have a common non-trivial irreducible. Then by taking a conjugate of
$H^{\prime}$ we can switch $\rho_{2,l}$ and $\rho_{2,l^{\prime}}$ so that the
$Y_{i,j}$ can be adjusted without affecting the last column of $Y.$ Now
suppose that even switching in not possible. Then it follows that $q$ is
trivial and that the number of orbits is greater than than $v.$ Define $X$ as
follows$.$%
\[
X=\left[
\begin{array}
[c]{ccccccc}%
E_{r_{1}} & 0 & \cdots & 0 & r_{1}y_{1}E_{r_{v+1}} & \cdots & r_{1}%
y_{l-v}E_{r_{l}}\\
0 & E_{r_{2}} & \cdots & 0 & r_{2}y_{1}E_{r_{v+1}} & \cdots & r_{2}%
y_{l-v}E_{r_{l}}\\
\vdots & \vdots & \ddots & 0 & \vdots & \ddots & \vdots\\
0 & 0 & 0 & E_{r_{v}} & r_{v}y_{1}E_{r_{v+1}} & \cdots & r_{v}y_{l-v}E_{r_{l}}%
\end{array}
\right]
\]
where $E_{n}$ is a row matrix of 1's, $r_{s}$ is the size of the $s$'th orbit
and
\begin{equation}
y_{1}r_{v+1}+\cdots+y_{l-v}r_{l}=-1,\text{ }y_{i}\neq0,\text{ }r_{i}\neq0.
\label{eq-sumriyi}%
\end{equation}
Assuming that the $y_{j}$'s satisfy the given constraints and at least one of
the $r_{i}$'s is non-zero then $X$ meets all the requirements. Since $p$ is
coprime to $\left\vert H\right\vert $, none of the $r_{i}$'s are zero. Set
$s=l-v\ >0$ and assume that $p\neq2.$ According to Proposition
\ref{pr-aiXieqzero} in the next section, only $\frac{(p-1)^{s}-(-1)^{s}}%
{p}+(-1)^{s}$ of the $(p-1)^{s}$ selections of $y_{i}$ sum to zero. Thus at
least one of the sums in equation \ref{eq-sumriyi} is non-zero and it may be
scaled to equal $-1$. If $p=2$ and $v>1$ then choose $X$ of the following
form
\[
X=\left[
\begin{array}
[c]{cccccccc}%
E_{r_{1}} & yE_{r_{1}} & \cdots & 0 & E_{r_{v+1}} & E_{r_{v+2}} & \cdots &
E_{r_{l}}\\
0 & E_{r_{1}} & \cdots & 0 & E_{r_{v+1}} & 0 & \cdots & 0\\
\vdots & \vdots & \ddots & 0 & \vdots & \vdots & \ddots & \vdots\\
0 & 0 & 0 & E_{r_{v}} & E_{r_{v+1}} & 0 & \cdots & 0
\end{array}
\right]  .
\]
Since all the $r_{i}$ are odd then $y$ may be chosen so that all row sums are
zero and that no column is zero.
\end{proof}

\subsection{Normalizers and conjugates\label{subsec-norm}}

Let us describe how to compute normalizers of subgroups fixing a point.
Suppose that $(g,\alpha)$ normalizes $H.$ Then for $(q(\beta),\beta)\in H,$
$(g,\alpha)(q(\beta),\beta)(g^{-1},\alpha^{-1})=(gq(\beta)g^{-1},\alpha
\beta\alpha^{-1})=(q(\alpha\beta\alpha^{-1}),\alpha\beta\alpha^{-1})\in H$ so
we must have
\begin{align}
\alpha\beta\alpha^{-1}  &  \in H^{\prime}\nonumber\\
gq(\beta)g^{-1}  &  =q(\alpha\beta\alpha^{-1}) \label{eq-norm}%
\end{align}
Thus we must already have $\alpha\in N_{\Sigma_{r}}(H^{\prime})=N^{\prime}.$
The normalizer $N^{\prime}$ permutes the representations of $H^{\prime}$ by
the formula
\[
q^{\alpha}(\beta)=q(\alpha\beta\alpha^{-1}),
\]
so there is $(g,\alpha)$ normalizing $H$ if and only if $q$ and
$q^{\alpha}$ are equivalent over $\mathbb{F}_{p}.$ In terms of characters,
$\alpha$ must fix the character of $q.$ Let $N^{\prime\prime}$ be the subgroup
satisfying $H^{\prime}\subseteq N^{\prime\prime}\subseteq N^{\prime}$ and
fixing the character of $q.$ If $(g_{1},\alpha),(g_{2},\alpha)$ both belong to
the normalizer of $H$ then equation (\ref{eq-norm}) implies $g_{1}%
q(\beta)g_{1}^{-1}=q(\alpha\beta\alpha^{-1})=g_{2}q(\beta)g^{-1}$ or
$q(\beta)=\left(  g_{1}^{-1}g_{2}\right)  q(\beta)\left(  g_{1}^{-1}%
g_{2}\right)  ^{-1}$ for all $\beta\in H^{\prime.}$ Thus,
$h=g_{1}^{-1}g_{2}\in Z=Z_{GL(v,\mathbb{F}_{p})}(q(H^{\prime})).$ It follows
that the normalizer of $H$ contains $Z$ as a subgroup with quotient
$N^{\prime\prime}$ and so the normalizer has size $\left\vert Z\right\vert
\left\vert N^{\prime\prime}\right\vert .$

Two subgroups fixing a point $H_{1}$ $H_{2}$ are conjugate if and only
$H_{2}^{\prime}=\alpha H_{2}^{\prime}a^{-1}$ and the two representations
$q^{\prime}$ and $q$ satisfy%
\[
q^{\prime}(\alpha\beta\alpha^{-1})=gq(\beta)g^{-1}%
\]
for some $g\in GL(v,\mathbb{F}_{p}).$

\subsection{Fixed point subsets and strata}

There are three ways we shall consider to calculate the size of fixed point subsets. We illustrate the first two by example with $v=2,$ $r=4$, $H^{\prime}$ =
$\left\langle (1,2)(3,4)\right\rangle $ and $q((1,2)(3,4))=\left[
\begin{array}
[c]{cc}%
1 & 0\\
0 & -1
\end{array}
\right]  .$ The first method is by brute force and can easily be  implemented by computer. 

\begin{example}
\label{ex-brute}
A typical $X=\left[
\begin{array}
[c]{cccc}%
x_{1} & x_{2} & x_{3} & x_{4}\\
y_{1} & y_{2} & y_{3} & y_{4}%
\end{array}
\right]  $ must satisfy the equations%
\begin{align*}
x_{1}+x_{2}+x_{3}+x_{4}=y_{1}+y_{2}+y_{3}+y_{4}=0  & \\
x_{1}-x_{2}=x_{3}+x_{4}=y_{1}-y_{2}=y_{3}+y_{4}=0  &
\end{align*}
and so $X=\left[
\begin{array}
[c]{cccc}%
x & x & -x & -x\\
y & -y & z & -z
\end{array}
\right]  .$ Clearly $x\neq0$ and at least one of the six $2\times2$ minors
must be non-zero, i.e., one of $xz,x(y-z),x(y+z),xy$ must be non- zero. It
follows that our conditions are $x\neq0$ and $(y,z)\neq(0,0)$ and the number
of valid $X$'s is $(p-1)(p^{2}-1).$

\end{example}

A drawback to this method is that the representation must be constructed. However, we note that it works if $p$ divides the order of a stabilizer.  The second method uses the "diagonalized" form of the intertwining operator
$Y=A^{-1}XB,$ and can usually be determined by inspection of the characters of
$q$ and $\overline{\rho}.$ 

\begin{example}
Observe that the character of $q$ is $\chi_{0}%
+\chi_{1}$ and $\overline{\rho}=\chi_{0}+2\chi_{1},$ where $\chi_{0}$ and
$\chi_{1}$ are the trivial and non-trivial characters of $H^{\prime}$
respectively. The matrices $A,B,$ and $Y$ can be chosen as:
\[
A=\left[
\begin{array}
[c]{cc}%
1 & 0\\
0 & 1
\end{array}
\right]  ,\text{ }B=\left[
\begin{array}
[c]{cccc}%
1 & 1 & 0 & 1\\
1 & -1 & 0 & 1\\
0 & 0 & 1 & 1\\
0 & 0 & -1 & 1
\end{array}
\right]  ,\text{ }Y=\left[
\begin{array}
[c]{cccc}%
x & 0 & 0 & 0\\
0 & y & z & 0
\end{array}
\right]  .
\]
We need $x\neq0$ and $(y,z)\neq(0,0)$ giving the same result as Example \ref{ex-brute}. Note
that we don't need to know $A$ and $B$ explicitly to determine $Y,$
just the characters.
\end{example}

The two previous methods are useful when the fixed point set has small
dimension and the non-vanishing polynomials are simple. At the other extreme
we need to calculate the fixed point set for the trivial point stabilizer. As a starting point we need to calculate the size of $\Omega.$ We can
calculate the number of vectors by inclusion-exclusion followed by a
specialized M\"{o}bius inversion. Let $V$ be an arbitrary vector space of
dimension $v$ over $\mathbb{F}_{p}$ and let
\begin{align*}
\overline{\Omega}^{r}(V)  &  =\left\{  X=\left(  X_{1},\ldots,X_{r}\right)
\in V^{r}:X_{i}\neq0,i=1,\ldots,r,\text{ }X_{1}+X_{2}+\cdots+X_{r}=0\right\}
\\
\Omega^{r}(V)  &  =\left\{  X\in\overline{\Omega}^{r}(V):\mathrm{\ }%
V=\left\langle X_{1},\ldots,X_{r}\right\rangle \right\}
\end{align*}
and set $\overline{\Omega}^{v,r}=\overline{\Omega}^{r}(\mathbb{F}_{p}^{v}),$
$\Omega^{v,r}=\Omega^{r}(\mathbb{F}_{p}^{v}).$ Clearly the cardinalities of
the sets depend only on $r$ and $v=\dim V.$

\begin{proposition}
Let $\overline{\Omega}^{r}(V)$ and $\Omega^{r}(V)$ be defined as above and
set
\[
\overline{\omega}(v,r)=\left\vert \overline{\Omega}^{r}(V)\right\vert ,\text{
}\omega(v,r)=\left\vert \Omega^{r}(V)\right\vert .
\]
Then,
\begin{align}
\overline{\omega}(v,r)  &  =\frac{(p^{v}-1)^{r}-(-1)^{r}}{p^{v}}%
+(-1)^{r}\label{eq-omegacount}\\
\omega(1,r)  &  =\overline{\omega}(1,r)\nonumber\\
\omega(v,r)  &  =\overline{\omega}(v,r)-\sum\limits_{l=1}^{v-1}n_{v,l}%
\omega(l,r)\nonumber
\end{align}
where
\[
n_{v,l}=\frac{\prod\limits_{j=l+1}^{v}\left(  p^{j}-1\right)  }{\prod
\limits_{j=1}^{v-l}\left(  p^{j}-1\right)  }%
\]
is the number of subspaces $L$ of dimension $l$ in $V.$
\end{proposition}

\begin{proof}
Let
\begin{align*}
\overline{\Omega}_{0}  &  =\left\{  \left(  X_{1},\ldots,X_{r}\right)  \in
V^{r}:\text{ }X_{1}+X_{2}+\cdots+X_{r}=0\right\} \\
\overline{\Omega}_{i}  &  =\left\{  \left(  X_{1},\ldots,X_{r}\right)  \in
V^{r}:X_{i}=0,\text{ }X_{1}+X_{2}+\cdots+X_{r}=0\right\}  ,i=1,\ldots,r\\
\overline{\Omega}  &  =\left\{  \left(  X_{1},\ldots,X_{r}\right)  \in
V^{r}:X_{i}\neq0,i=1,\ldots,r,\text{ }X_{1}+X_{2}+\cdots+X_{r}=0\right\}
\end{align*}
so that $\overline{\Omega}^{r}(V)=\overline{\Omega}=\overline{\Omega}%
_{0}-\bigcup\limits_{i}\overline{\Omega}_{i}.$ Then, by inclusion-exclusion:%
\begin{align*}
\overline{\omega}(v,r)  &  =\left\vert \overline{\Omega}\right\vert \\
&  =\left\vert \overline{\Omega}_{0}\right\vert -\sum_{i}\left\vert
\overline{\Omega}_{i}\right\vert +\sum_{i<j}\left\vert \overline{\Omega}%
_{i}\cap\overline{\Omega}_{j}\right\vert -\sum_{i<j<k}\left\vert
\overline{\Omega}_{i}\cap\overline{\Omega}_{j}\cap\overline{\Omega}%
_{k}\right\vert +\cdots\\
&  +(-1)^{r}\left\vert \overline{\Omega}_{i}\cap\cdots\cap\overline{\Omega
}_{r}\right\vert \\
&  =p^{v(r-1)}-\binom{r}{1}p^{v(r-2)}+\binom{r}{2}p^{v(r-3)}-\binom{r}%
{3}p^{v(r-3)}+\cdots\\
&  +(-1)^{r-1}\binom{r}{r-1}+(-1)^{r}\\
&  =\frac{(p^{v}-1)^{r}-(-1)^{r}}{p^{v}}+(-1)^{r}%
\end{align*}
Next let $L$ be denote an arbitrary proper subspace of $V$ of dimension $l.$
Then%
\[
\overline{\Omega}^{r}(V)=\Omega^{r}(V)\cup\bigcup\limits_{0\subset L\subset
V}\Omega^{r}(L)
\]
It follows that
\begin{align*}
\overline{\omega}(v,r)  &  =\left\vert \overline{\Omega}^{r}(V)\right\vert
=\left\vert \Omega^{r}(V)\right\vert +\sum_{l=1}^{v-1}\sum\limits_{\dim
L=l}\Omega^{r}(L)\\
&  =\left\vert \Omega^{r}(\mathbb{F}_{p}^{v})\right\vert +\sum_{l=1}%
^{v-1}n_{v,l}\left\vert \Omega^{r}(\mathbb{F}_{p}^{l})\right\vert \\
&  =\omega(v,r)+\sum_{l=1}^{v-1}n_{v,l}\omega(l,r)
\end{align*}
The remaining formulas in equations (\ref{eq-omegacount}) follow immediately.
Finally, by noting that $\left\vert GL_{v}(\mathbb{F}_{p})\right\vert
=\prod\limits_{j=1}^{v}\left(  p^{j}-1\right)  p^{(n^{2}-n)/2}$, by a
homogeneous space argument, letting $GL_{v}(\mathbb{F}_{p})$ act on the
subspaces, we have 
\[
n_{v,l}=\frac{\left\vert GL_{v}(\mathbb{F}_{p})\right\vert }{\left\vert
GL_{l}(\mathbb{F}_{p})\right\vert \left\vert GL_{v-l}(\mathbb{F}%
_{p})\right\vert p^{l(v-l)}}=\frac{\prod\limits_{j=l+1}^{v}\left(
p^{j}-1\right)  }{\prod\limits_{j=1}^{v-l}\left(  p^{j}-1\right)  }.
\]

\end{proof}

\begin{example}
For later work we will need the following values $\left\vert \Omega
^{r}(\mathbb{F}_{p}^{v})\right\vert .$
\[%
\begin{tabular}
[c]{|c|c|c|}\hline
$v$ & $r$ & $\left\vert \Omega^{r}(\mathbb{F}_{p}^{v})\right\vert $\\\hline
$1$ & $2$ & $p-1$\\\hline
$1$ & $3$ & $\left(  p-1\right)  \left(  p-2\right)  $\\\hline
$1$ & $4$ & $\left(  p-1\right)  \left(  p^{2}-3p+3\right)  $\\\hline
$1$ & $5$ & $\left(  p-1\right)  (p^{3}-4p^{2}+6p-4)$\\\hline
$2$ & $3$ & $p\left(  p-1\right)  \left(  p^{2}-1\right)  $\\\hline
$2$ & $4$ & $p\left(  p-1\right)  \left(  p^{2}-1\right)  \left(
p^{2}+p-3\right)  $\\\hline
$2$ & $5$ & $p\left(  p-1\right)  \left(  p^{2}-1\right)  (p^{4}+p^{3}%
-3p^{2}-4p+6)$\\\hline
$3$ & $4$ & $p^{3}\left(  p-1\right)  \left(  p^{2}-1\right)  \left(
p^{3}-1\right)  $\\\hline
\end{tabular}
\
\]
Note that, since $GL_{v}(\mathbb{F}_{p})$ acts freely by left multiplication
on $\Omega^{r}(\mathbb{F}_{p}^{v})$ then $\left\vert \Omega^{r}(\mathbb{F}%
_{p}^{v})\right\vert $ is always evenly divided by $\left\vert GL_{v}%
(\mathbb{F}_{p})\right\vert .$ The $\left\vert \Omega^{r}(\mathbb{F}_{p}%
^{v})\right\vert $ have been written to show this fact (some can be further factored).
\end{example}

\begin{example}
\label{ex-strat1}Assume that $q:H^{\prime}\rightarrow GL_{v}(\mathbb{F}_{p}) $
is trivial. Then $H$ orbits determine a partition of $r$, $r=r_{1}%
+r_{2}+\cdots r_{s}$ with $1\leq r_{1}\leq r_{2}\leq$ $\cdots\leq r_{s}.$ Such
a partition of $r$ determines $\binom{r}{r_{1},\ldots,r_{s}}$ partitions of
$\{1,\ldots r\}.$ with the canonical one being
\[
\left\{  \{1,\ldots,r_{1}\},\{r_{1}+1,\ldots,r_{1}+r_{2}\},\ldots
,\{r_{1}+\cdots r_{s-1}+1,\ldots,r\}\right\}
\]
A typical element of the stratum is obtained by selecting distinct non-zero
$X_{1},\ldots X_{s}$ satisfying
\begin{align*}
r_{1}X_{1}+\cdots+r_{s}X_{s}  &  =0\\
\mathrm{rank}\left[
\begin{array}
[c]{ccc}%
X_{1} & \cdots & X_{s}%
\end{array}
\right]   &  =v
\end{align*}
and then taking in order $r_{1}$ of $X_{1},r_{2}$ of $X_{2},\ldots r_{s}$ of
$X_{s}.$ To count the points fixed by $H$ we need the following generalization
\end{example}

\begin{proposition}
\label{pr-aiXieqzero}Suppose that $a_{1},\ldots,a_{s}$ are non-zero scalars in
$\mathbb{F}_{p}.$ Then the cardinality of the set
\[
\left\{
\begin{array}
[c]{c}%
X=\left[
\begin{array}
[c]{ccc}%
X_{1} & \cdots & X_{s}%
\end{array}
\right]  \in\mathbb{F}_{p}^{vs}:X_{i}\neq0,i=1,\ldots,s,\text{ }\\
a_{1}X_{1}+\cdots+a_{s}X_{s}=0,\text{ }\mathrm{rank}X=v
\end{array}
\right\}
\]
is $\left\vert \Omega^{v,s}\right\vert .$
\end{proposition}

\begin{proof}
The map $\left[
\begin{array}
[c]{ccc}%
X_{1} & \cdots & X_{s}%
\end{array}
\right]  \longleftrightarrow\left[
\begin{array}
[c]{ccc}%
a_{1}X_{1} & \cdots & a_{s}X_{s}%
\end{array}
\right]  $ is a bijection between the given set and $\Omega^{v,r}$.
\end{proof}

\subsection{Singular Primes and Galois fusion}

The forgoing representation theory depends on certain divisibility properties
of primes. The first problem occurs when a prime $p$ divides the order of
$H^{\prime}.$ Let us call such a prime a singular prime. The singular primes
must satisfy $p\leq r$. The main problem is that complete reducibility of the
representations $q$ and $\rho$ fails. Moreover there is new twist in that the
sum over an orbit may be zero if the size of the orbit is divisible by the
prime. In this case one needs to directly compute the representations and then compute the normalizers and fixed point set by brute force.

Another problem that can occur is that $\mathbb{F}_{p}$ may not be a splitting
field for the representations $H^{\prime}.$ In trying to describe the
representations completely in terms of characters will have to be take
characters together in Galois equivalent groups in order that the
representations $q$ be defined over the base field $\mathbb{F}_{p}$. We will
call this \emph{Galois fusion}. Furthermore, the formulas for the number of
fixed points will have different polynomials in the primes. To help in
describing this theory we use the following notation. For a prime $p$ and an
integer $n$ let $U_{n}(\mathbb{F}_{p})$ be the set of primitive $n^{\prime}$th
roots of 1 in $\mathbb{F}_{p}$ and let $\beta_{n}(p)=\left\vert U_{n}%
(\mathbb{F}_{p})\right\vert .$ Again, instead of developing the full theory
here we will illustrate the ideas with the sample calculations in the next subsection.

\subsection{ Rank 2 actions with 4 branch points}

Our general approach to compute the number of orbits with given $p$-rank and
number of branch points is as follows.

\begin{enumerate}
\item Determine all subgroups fixing at least one point, one representative
for each conjugacy class. List them as $H_{1},H_{2},\ldots,H_{s}$ as described
in Subsection \ref{subsec-matrix}. The sequence of groups will depend upon the prime.

\item Compute $\left\vert \Omega^{H}\right\vert $ and $\left\vert
N_{G}(H)\right\vert $ for each subgroup.

\item Compute the matrix $D=\left[  d_{i,j}\right]  $ for the sequence of groups.

\item Compute $\left\vert \left\langle \mathring{\Omega}^{H}\right\rangle
/G\right\vert $ for each subgroup using formula $\ref{eq-orbcount}$ and the
add up the results.
\end{enumerate}

We will illustrate the steps by giving complete details for $p$-rank 2 and 4
branch points. For Step 1 we use a computer algebra system such as Magma or
GAP to do the following:

\begin{itemize}
\item Determine a representative $H^{^{\prime}}$of each conjugacy class of
$\Sigma_{4}.$

\item Determine the character table of $H^{\prime}.$

\item Decompose the reduced natural representation of $H^{\prime}$ into
irreducibles over $\mathbb{C}.$ Determine all compatible representations of
$H^{\prime}$ and how they reduce over $\mathbb{F}_{p}.$
\end{itemize}

\begin{remark}
For each modular representation of an arbitrary group $G$ defined over $p$
there is complex valued Brauer character defined on the $p$-regular elements
of $G$ (order coprime to $p$). In the case $p\nmid\left\vert G\right\vert $ each
Brauer character is an ordinary character and the modular representations of
$G$ can be completely described by the irreducible characters. Thus in the
discussion below we describe the modular characters by sums of irreducible
ordinary characters from the character table of $G.$
\end{remark}

\bigskip

\noindent\textbf{Step 1. }There are 11 conjugacy classes of subgroups of
$\Sigma_{4}.$ Using Magma, the orbits, and the decomposition of the reduced
permutation (in terms of the characters of the subgroup), and the degrees of
the representation are computed. The column giving the reduced permutation
representation is written in terms of the characters of the subgroup using the
order in the character table produced by Magma. The list of the degrees is in
the order given by the character table. $\ $The following notation is used:
$\Sigma_{k},A_{k}$  denote respectively the symmetric and alternating groups on $k$ points, $C_{k_{1},\ldots,k_{s}}$ the cyclic group with cycle structure $(k_{1}%
,\ldots,k_{s})$, $D_{k}$ dihedral group on $k$ points and $V_{4}$ denotes the Klein 4 group.%

\[%
\begin{tabular}
[c]{c}%
Table 4.1 Subgroups of $\Sigma_{4}$\\
$%
\begin{tabular}
[c]{|c|c|c|c|c|}\hline
$H^{\prime}$ & $\left\vert H^{\prime}\right\vert $ & orbits & red perm rep &
degrees\\\hline
$\Sigma_{1}^{4}$ & $1$ & $\{1\},\{2\},\{3\},\{4\}$ & $3\chi_{1}$ & $1$\\\hline
$\Sigma_{2}\times\Sigma_{1}^{2}$ & $2$ & $\{1,2\},\{3\},\{4\}$ & $2\chi
_{1}+\chi_{2}$ & $1,1$\\\hline
$C_{2,2}$ & $2$ & $\{1,2\},\{3,4\}$ & $\chi_{1}+2\chi_{2}$ & $1,1$\\\hline
$A_{3}\times\Sigma_{1}$ & $3$ & $\{1,2,3\},\{4\}$ & $\chi_{1}+\chi_{2}%
+\chi_{3}$ & $1,1,1$\\\hline
$\Sigma_{2}\times\Sigma_{2}$ & $4$ & $\{1,2\},\{3,4\}$ & $\chi_{1}+\chi
_{3}+\chi_{4}$ & $1,1,1,1$\\\hline
$V_{4}$ & $4$ & $\{1,2,3,4\}$ & $\chi_{2}+\chi_{3}+\chi_{4}$ & $1,1,1,1$%
\\\hline
$C_{4}$ & $4$ & $\{1,2,3,4\}$ & $\chi_{2}+\chi_{3}+\chi_{4}$ & $1,1,1,1$%
\\\hline
$\Sigma_{3}\times\Sigma_{1}$ & $6$ & $\{1,2,3\},\{4\}$ & $\chi_{1}+\chi_{3}$ &
$1,1,2\ $\\\hline
$D_{4}$ & $8$ & $\{1,2,3,4\}$ & $\chi_{3}+\chi_{5}$ & $1,1,1,1,2$\\\hline
$A_{4}$ & $12$ & $\{1,2,3,4\}$ & $\chi_{4}$ & $1,1,1,3$\\\hline
$\Sigma_{4}$ & $24$ & $\{1,2,3,4\}$ & $\chi_{4}$ & $1,1,2,3,3$\\\hline
\end{tabular}
\ \ $%
\end{tabular}
\ \
\]

Next we determine representatives of conjugacy classes of subgroups fixing a
point, given in Table 4.2. Using Proposition \ref{pr-stabchar}, a stabilizer
can be specified by a compatible character if one exists. Note that here we
need to take away a one-dimensional character from the reduced representation
to form $q.\ $This automatically eliminates $\Sigma_{4},$ $A_{4}$ and
$\Sigma_{3}\times\Sigma_{1}.$ If there is a fixed point of $H^{\prime}$ then
$q$ must contain the trivial representation $\chi_{1}.$ This eliminates
$\Sigma_{3}\times\Sigma_{1}$ and taking $\chi_{1}$ away from the reduced
representation of $A_{3}\times\Sigma_{1}.$ Also in the $A_{3}\times\Sigma_{1}$
case, since we can take only one $\chi_{2}$ or $\chi_{3}$ then $\mathbb{F}%
_{p}$ must contain primitive cube roots of 1. Different characters may lead to
conjugate stabilizers. By the discussion in Subsection \ref{subsec-norm} if
$\alpha$ normalizes $H^{\prime},$ and $g\in$ $GL(2,\mathbb{F}_{p})$ then the
subgroup $H^{\prime}$ and $gq^{\alpha}g^{-1}$ determines a conjugate subgroup.
The $N_{H^{\prime}}$\emph{ }$\chi$\emph{-orbits}\ column in Table 4.2 lists
the non-trivial orbits of characters under conjugation. If two linear
combinations of characters are equivalent under the normalizer then conjugate
stabilizers are determined and so we need only write one down. Furthermore,
depending on the number of roots of unity in $\mathbb{F}_{p}$ some
representations may not reduce completely. In that case the characters have to
be taken together as a group, namely the set of all Galois conjugate
characters. The non-trivial Galois orbits are noted in the \emph{Galois
fusion}\ column. Cases 11 and 11a are the same subgroup but are distinguished
because $\left\vert \Omega^{H_{i}}\right\vert ,$ Table 4.3, depends on the
value of $\beta_{4}(p).$%

\[%
\begin{tabular}
[c]{c}%
Table 4.2 subgroups fixing a point, $p\neq2,3$\\
$%
\begin{tabular}
[c]{|c|c|c|c|c|c|}\hline
Case & $H^{\prime}$ & $q$ & $N_{H^{\prime}}$ $\chi$-orbits & Galois Fusion &
restrictions\\\hline
$1$ & $\Sigma_{1}^{4}$ & $2\chi_{1}$ &  &  & \\\hline
$2$ & $\Sigma_{2}\times\Sigma_{1}^{2}$ & $2\chi_{1}$ &  &  & \\\hline
$3$ & $\Sigma_{2}\times\Sigma_{1}^{2}$ & $\chi_{1}+\chi_{2}$ &  &  & \\\hline
$4$ & $C_{2,2}$ & $\chi_{1}+\chi_{2}$ &  &  & \\\hline
$5$ & $C_{2,2}$ & $2\chi_{2}$ &  &  & \\\hline
$6$ & $A_{3}\times\Sigma_{1}$ & $\chi_{1}+\chi_{2}$ & $\left\{  \chi_{2}%
,\chi_{3}\right\}  $ &  & $\beta_{3}(p)=2$\\\hline
$7$ & $\Sigma_{2}\times\Sigma_{2}$ & $\chi_{1}+\chi_{3}$ & $\left\{  \chi
_{3},\chi_{4}\right\}  $ &  & \\\hline
$8$ & $\Sigma_{2}\times\Sigma_{2}$ & $\chi_{3}+\chi_{4}$ & $\left\{  \chi
_{3},\chi_{4}\right\}  $ &  & \\\hline
$9$ & $K_{4}$ & $\chi_{2}+\chi_{3}$ & $\left\{  \chi_{2},\chi_{3},\chi
_{4}\right\}  $ &  & \\\hline
$10$ & $C_{4}$ & $\chi_{2}+\chi_{3}$ & $\left\{  \chi_{2},\chi_{4}\right\}  $
&  & $\beta_{4}(p)=2$\\\hline
$11$ & $C_{4}$ & $\chi_{2}+\chi_{4}$ & $\left\{  \chi_{2},\chi_{4}\right\}  $
& $\left\{  \chi_{2},\chi_{4}\right\}  $ & $\beta_{4}(p)=0$\\\hline
$11a$ & $C_{4}$ & $\chi_{2}+\chi_{4}$ & $\left\{  \chi_{2},\chi_{4}\right\}  $
&  & $\beta_{4}(p)=2$\\\hline
$13$ & $D_{4}$ & $\chi_{5}$ &  &  & \\\hline
\end{tabular}
\ \ $%
\end{tabular}
\ \
\]

\bigskip

\noindent\textbf{Step 2.} Next we list $\left\vert N_{G}(H_{i})\right\vert $
and the $\left\vert \Omega^{H_{i}}\right\vert .$ Each calculation of
$\left\vert \Omega^{H_{i}}\right\vert $ was confirmed by using the brute force
method using Maple. We have the following table:%

\[%
\begin{tabular}
[c]{c}%
Table 4.3 fixed point data for $\Sigma_{4},$ $p\neq2,3$\\
$%
\begin{tabular}
[c]{|c|c|c|c|}\hline
Case & $\left\vert H^{\prime}\right\vert $ & $N_{G}(H_{i})$ & $\left\vert
\Omega^{H_{i}}\right\vert $\\\hline
$1$ & $1$ & $24p\left(  p-1\right)  \left(  p^{2}-1\right)  $ & $p\left(
p-1\right)  \left(  p^{2}-1\right)  \left(  p^{2}+p-3\right)  $\\\hline
$2$ & $2$ & $4p\left(  p-1\right)  \left(  p^{2}-1\right)  $ & $p\left(
p-1\right)  \left(  p^{2}-1\right)  $\\\hline
$3$ & $2$ & $4(p-1)^{2}$ & $(p-1)\left(  p^{2}-1\right)  $\\\hline
$4$ & $2$ & $8(p-1)^{2}$ & $(p-1)\left(  p^{2}-1\right)  $\\\hline
$5$ & $2$ & $8p\left(  p-1\right)  \left(  p^{2}-1\right)  $ & $p\left(
p-1\right)  \left(  p^{2}-1\right)  $\\\hline
$6$ & $3$ & $3(p-1)^{2}$ & $\left(  p-1\right)  ^{2}$\\\hline
$7$ & $4$ & $4(p-1)^{2}$ & $\left(  p-1\right)  ^{2}$\\\hline
$8$ & $4$ & $8(p-1)^{2}$ & $\left(  p-1\right)  ^{2}$\\\hline
$9$ & $4$ & $8(p-1)^{2}$ & $\left(  p-1\right)  ^{2}$\\\hline
$10$ & $4$ & $4(p-1)^{2}$ & $\left(  p-1\right)  ^{2}$\\\hline
$11$ & $4$ & $8(p-1)^{2}$ & $p^{2}-1$\\\hline
$11a$ & $4$ & $8(p^{2}-1)$ & $\left(  p-1\right)  ^{2}$\\\hline
$12$ & $10$ & $8(p-1)$ & $p-1$\\\hline
\end{tabular}
\ \ \ $%
\end{tabular}
\ \ \
\]

\bigskip

\noindent\textbf{Step 3.} To compute the $d_{i,j}$'s, let us first look at the
case where all roots of unity are available $\beta_{3}(p)=\beta_{4}%
(p)=2.\ $The other cases may be derived from this case. We present the
$d_{i,j}$'s, in matrix form with this ordering of the groups.%

\[%
\begin{tabular}
[c]{c}%
Table 4.4 Subgroup ordering\\
$%
\begin{tabular}
[c]{|c|c|c|c|c|c|c|}\hline
$i$ & $1$ & $2$ & $3$ & $4$ & $5$ & $6$\\\hline
Case & $1$ & $2$ & $3$ & $4$ & $5$ & $6$\\\hline
$H_{i}^{\prime}$ & $\Sigma_{1}^{4}$ & $\Sigma_{2}\times\Sigma_{1}^{2}$ &
$\Sigma_{2}\times\Sigma_{1}^{2}$ & $C_{2,2}$ & $C_{2,2}$ & $A_{3}\times
\Sigma_{1}$\\\hline
$q_{i}$ & $2\chi_{1}$ & $2\chi_{1}$ & $\chi_{1}+\chi_{2}$ & $\chi_{1}+\chi
_{2}$ & $2\chi_{2}$ & $\chi_{1}+\chi_{2}$\\\hline
\end{tabular}
\ \ \ $\\
\\
$%
\begin{tabular}
[c]{|c|c|c|c|c|c|c|}\hline
$i$ & $7$ & $8$ & $9$ & $10$ & $11$ & $12$\\\hline
Case & $7$ & $8$ & $9$ & $10$ & $11a$ & $12$\\\hline
$H_{i}^{\prime}$ & $\Sigma_{2}\times\Sigma_{2}$ & $\Sigma_{2}\times\Sigma_{2}$
& $K_{4}$ & $C_{4}$ & $C_{4}$ & $D_{4}$\\\hline
$q_{i}$ & $\chi_{1}+\chi_{3}$ & $\chi_{3}+\chi_{4}$ & $\chi_{2}+\chi_{3}$ &
$\chi_{2}+\chi_{3}$ & $\chi_{3}+\chi_{4}$ & $\chi_{5}$\\\hline
\end{tabular}
\ \ \ $%
\end{tabular}
\ \ \ \
\]
The values of $d_{i,j}=\left\vert D(H_{i},H_{j})\right\vert $ can be computed
as follows$.$ Using Magma or GAP, find all the subgroup classes of
$H_{j}^{\prime}.$ For each subgroup class $\left\langle K\right\rangle ,$
$K\subseteq H_{j}^{\prime}$ find the corresponding $H_{i}$'s such that
$H_{i}^{\prime}$ is conjugate to $K.$ Transport the character $q_{j}$ of
$H_{j}$ to $H_{i}^{\prime}$ to determine the subgroup $H_{i}.$ The entry for
$d_{i,j}$ is then $\left\vert \left\langle K\right\rangle \right\vert $
(computed in $H_{j}^{\prime}$)$.$

\begin{example}
Here is the matrix $D$ for the set of subgroups above.%
\[
\left[
\begin{array}
[c]{ccccccccccccc}
& H_{1} & H_{2} & H_{3} & H_{4} & H_{5} & H_{6} & H_{7} & H_{8} & H_{9} &
H_{10} & H_{11} & H_{12}\\
H_{1} & 1 & 1 & 1 & 1 & 1 & 1 & 1 & 1 & 1 & 1 & 1 & 1\\
H_{2} &  & 1 & 0 & 0 & 0 & 0 & 1 & 0 & 0 & 0 & 0 & 0\\
H_{3} &  &  & 1 & 0 & 0 & 0 & 1 & 2 & 0 & 0 & 0 & 2\\
H_{4} &  &  &  & 1 & 0 & 0 & 1 & 0 & 2 & 1 & 0 & 2\\
H_{5} &  &  &  &  & 1 & 0 & 0 & 1 & 1 & 0 & 1 & 1\\
H_{6} &  &  &  &  &  & 1 & 0 & 0 & 0 & 0 & 0 & 0\\
H_{7} &  &  &  &  &  &  & 1 & 0 & 0 & 0 & 0 & 0\\
H_{8} &  &  &  &  &  &  &  & 1 & 0 & 0 & 0 & 1\\
H_{9} &  &  &  &  &  &  &  &  & 1 & 0 & 0 & 1\\
H_{10} &  &  &  &  &  &  &  &  &  & 1 & 0 & 0\\
H_{11} &  &  &  &  &  &  &  &  &  &  & 1 & 1\\
H_{12} &  &  &  &  &  &  &  &  &  &  &  & 1
\end{array}
\right]
\]

\end{example}

\bigskip

\noindent\textbf{Step 4.} Finally we compute the number of orbits according to
formula (\ref{eq-orbcount}). We have $O^{\circ}=TD^{-1}N^{-1}L$ where$\ T$ and
$N$ are the diagonal matrices formed from columns 2 and 3 of Table 4.3,
respectively, and $L$ is the vector formed from column 4 of Table 4.3.
Accordingly we obtain the vector of orbits and the total number of orbits%
\begin{align*}
O^{\circ}  &  =\left[
\begin{array}
[c]{cccccccccccc}%
\frac{p^{2}-8p+7}{24} & 0 & \frac{p-3}{2} & \frac{p-5}{4} & 0 & 1 & 1 & 0 &
0 & 1 & 0 & 1
\end{array}
\right]  ^{\top}\\
\left\vert \Omega/G\right\vert  &  =\frac{p^{2}+10p+37}{24}%
\end{align*}
We may check that all the numbers are integers by substituting in $p=12k+1,$
since $12$ divides $p-1.$

There are a total of 4 families of primes which depend upon whether $\beta
_{4}(p)=0,2$ and $\beta_{3}(p)=0,2.$ For $p>3$ these families depend on the
congruence class of $p\operatorname{mod}12.$ For each family certain subgroup
cases from Table 4.1 are excluded. The corresponding matrices $D$ are obtained
from the sample matrix above by deleting the row(s) and column(s)
corresponding to the excluded subgroups Depending on the whether $\beta
_{4}(p)=0,2$ one chooses the same subgroup $H_{10}$ but with a different
normalizer and fixed point data given by cases 11 and 11a.
\[%
\begin{tabular}
[c]{|c|c|c|c|c|}\hline
$\beta_{3}(p)$ & $\beta_{4}(p)$ & $p\operatorname{mod}12$ & excluded
subgroups & $\left\vert \Omega/G\right\vert $\\\hline
$0$ & $0$ & $11$ & $H_{6},H_{10}$ & $\frac{1}{24}(p^{2}+6p+9)$\\\hline
$2$ & $0$ & $7$ & $H_{10}$ & $\frac{1}{24}(p^{2}+6p+25)$\\\hline
$0$ & $2$ & $5$ & $H_{6}$ & $\frac{1}{24}(p^{2}+6p+31)$\\\hline
$2$ & $2$ & $1$ & none & $\frac{1}{24}(p^{2}+6p+37)$\\\hline
\end{tabular}
\ \ \
\]
The complete results for 3 and 4 branch points are given in Table 5.1.

\section{Low genus examples\label{Sec-Examples}}

To finish, we illustrate our results through a number of interesting examples.
To emphasize the explicit results which can be obtained, the first example we
consider is very specific. Following this, we shall present some general
examples which hold independent of the prime $p$.

\begin{example}
\label{ex-genus26}Let $G$ be an elementary abelian group of order $25$. Using
the Riemann-Hurwitz formula, it can be shown that there are exactly two
different signatures for $\Gamma$ for a surface of genus $26$ which give rise
to subgroups of $\mathcal{M}_{26}$ isomorphic to $G$ - the signatures $(2,-)$
and $(1;5,5,5)$. To find the total number of conjugacy classes of subgroups of
$\mathcal{M}_{26}$ isomorphic to $G$, we need to find the number of classes
induced by each of these signatures.

For the signature $(2,-)$, a direct application of Corollary \ref{cor-Number}
gives two different classes of subgroups. For the signature $(1;5,5,5)$, there
are two possibilities for the group $G^{h}$, either trivial or cyclic of order
$5$. If $G^{h}$ is trivial, $G^{e}=G$, so Proposition \ref{pr-onemore} implies
there is just one epimorphism from $\Gamma$ onto $G$. Else, the image is
cyclic of order $5$. In this case, we apply Corollary \ref{cor-counting}. For
the hyperbolic part, observe that Corollary \ref{cor-Number} tells us there
will be a single epimorphism arising from $\Gamma$ with signature $(1;-)$ onto
$G^{h}$. For the elliptic part, Example \ref{ex-Small} below shows that all
epimorphisms from $\Gamma$ with signature $(0;5,5,5)$ onto $G^{e}$ are
equivalent. Thus there is just one epimorphism arising from the elliptic part.
Therefore, in total, there are $2+1+1\ast1=4$ conjugacy classes of subgroups
isomorphic to $G$ in $\mathcal{M}_{26}$.
\end{example}

\begin{example}
\label{ex-Small} We can use our results to enumerate the equivalence classes of totally ramified actions for 3 or 4 branch points. Complete results for these cases are given in the following table: 
\[%
\begin{tabular}
[c]{c}%
Table 5.1 $r=$ \# branch points, $v=\ p$-rank\\
$%
\begin{tabular}
[c]{|c|c|c|c|c|}\hline
$r$ & $v$ & primes & congruence & \#equiv classes\\\hline
$3$ & $1$ & $2$ &  & $0$\\\hline
$3$ & $1$ & $3$ &  & $1$\\\hline
$3$ & $1$ & $\beta_{3}(p)=0$ & $p=5\operatorname{mod}6$ & $\frac{1}{6}\left(
p+1\right)  $\\\hline
$3$ & $1$ & $\beta_{3}(p)=2$ & $p=1\operatorname{mod}6$ & $\frac{1}{6}\left(
p+5\right)  $\\\hline
$3$ & $2$ & $\mathrm{all}$ &  & $1$\\\hline
$4$ & $1$ & $2,3$ &  & $1$\\\hline
$4$ & $1$ & $\beta_{4}(p)=0$ & $p=7,11\operatorname{mod}12$ & $\frac{1}%
{24}(p^{2}+6p+5)$\\\hline
$4$ & $1$ & $\beta_{4}(p)=2$ & $p=1,5\operatorname{mod}12$ & $\frac{1}%
{24}(p^{2}+6p+17)$\\\hline
$4$ & $2$ & $3$ &  & $2$\\\hline
$4$ & $2$ & $\beta_{3}(p)=0,\beta_{4}(p)=0$ & $p=11\operatorname{mod}12$ &
$\frac{1}{24}(p^{2}+10p+9)$\\\hline
$4$ & $2$ & $\beta_{3}(p)=2,\beta_{4}(p)=0$ & $p=7\operatorname{mod}12$ &
$\frac{1}{24}(p^{2}+10p+25)$\\\hline
$4$ & $2$ & $\beta_{3}(p)=0,\beta_{4}(p)=2$ & $p=5\operatorname{mod}12$ &
$\frac{1}{24}(p^{2}+10p+21)$\\\hline
$4$ & $2$ & $\beta_{3}(p)=2,\beta_{4}(p)=2$ & $p=1\operatorname{mod}12$ &
$\frac{1}{24}(p^{2}+10p+37)$\\\hline
$4$ & $3$ & $\mathrm{all}$ &  & $1$\\\hline
\end{tabular}
\ \ \ \ \ \ $%
\end{tabular}
\ \ \ \ \ \
\]

\end{example}

The techniques and results developed in the previous examples allow us to
develop some further results for general families of signatures and groups.

\begin{example}
Suppose $\Gamma$ has signature $(g_{1};p,p,p)$ with $p>3$ and let $G$ denote
an elementary abelian group of order $p^{2}$. In this case, $G^{h}$ is either
trivial or has order $p$. We first consider the case when $G^{h}$ has order
$p$. By Example \ref{ex-Small}, there will be $(p+1+2\beta_{3}(p))/6$ totally
ramified epimorphisms onto $C_{p}$ and Corollary \ref{cor-Number} implies
there is a unique hyperbolic epimorphism onto $C_{p}$ giving a total of
$1+(p+1+2\beta_{3}(p))/6$. If $G^{h}$ is trivial, then Proposition
\ref{pr-onemore} implies there exists a unique elliptic epimorphism. Thus
there is a total of
\[
\frac{p+7+2\beta_{3}(p)}{6}%
\]
conjugacy classes of groups induced by this signature.
\end{example}

The results we have obtained for elementary abelian subgroups can in certain special cases be extended to provide information about other groups. We illustrate with the following two examples.

\begin{example}
By Proposition \ref{pr-onemore}, there is a unique conjugacy class of
elementary abelian subgroups of order $p^{2}$ in $\mathcal{M}_{(p-1)(p-2)/2}$
 with signature $(0;p,p,p)$. In fact, this can be
derived explicitly by showing that any two epimorphisms $\eta_{1},\eta
_{2}\colon\Gamma\rightarrow G$ from $\Gamma$ with signature $(0;p,p,p)$ differ by an automorphism $\alpha\in\aut{(G)}$,
$\eta_{1}=\alpha\circ\eta_{2}$. In \cite{Woo}, it is shown that the kernel
$\Ker{(\eta_{1} )}$ is in fact a uniformizing surface group for the $p$th Fermat
curve with defining equation $x^{p}+y^{p}=1$. The full automorphism group of
the $p$th Fermat curve is isomorphic to the semi-direct product $S_{3}%
\ltimes(\mathbb{F}_{p}\times\mathbb{F}_{p})\cong\Gamma_{1}/\Ker{(\eta_{1})}$ where
$\Gamma_{1}$ has signature $(0;2,3,2p)$. By the uniqueness of $\Ker{(\eta_{1})}$
and the uniqueness of the elementary abelian subgroup of $S_{3}\ltimes
(\mathbb{F}_{p}\times\mathbb{F}_{p})$, it follows that any two epimorphisms
$\eta_{1},\eta_{2}\colon\Gamma_{1}\rightarrow S_{3}\ltimes(\mathbb{F}%
_{p}\times\mathbb{F}_{p})$ must differ by an automorphism $\alpha
\in\aut{(S_{3}
\ltimes (\mathbb{F}_{p}\times \mathbb{F}_{p}))}$ . In particular, there will
be a unique conjugacy class of subgroups of $\mathcal{M}_{(p-1)(p-2)/2}$
isomorphic to $S_{3}\ltimes(\mathbb{F}_{p}\times\mathbb{F}_{p})$ induced by
$\Gamma$ with signature $(0;2,3,2p)$.
\end{example}

\begin{example}
\label{ex-acc} Let $G$ be an elementary abelian group of order $p^{2}$ where
$p\geq5$ and suppose that $\Gamma$ has signature $(0;p,p,p,p)$. In \cite{Woo}
it is shown that there exists a unique epimorphism $\eta_{1}$ from $\Gamma
_{1}$ with signature $(0;2,2p,4)$ onto $H=D_{4}\ltimes(C_{p}\times C_{p})$ (up
to the action of $\aut{(H)}$) which restricts to an epimorphism $\eta_{1}|_{\Gamma}\colon\Gamma\rightarrow C_{p}\times C_{p}$.
 In particular, this implies there exists a unique class of subgroups of $\mathcal{M}_{(p-1)^2}$ isomorphic to $H$ with signature $(0;2,2p,4)$.  Alternatively, this can also be seen on the level of generating
vectors of the restriction $\eta_{1}|_{\Gamma}$. Specifically, it can be shown that if $x$
and $y$ generate $G$, then all generating vectors for epimorphisms which extend to epimorphisms from $\Gamma_{1}$ to $H$ are  $\aut{(G)}\times
\aut{(\Gamma )}$ equivalent to $(x,x^{-1},y,y^{-1})$. This implies there exists a unique class of groups isomorphic to $G$ in $\mathcal{M}_{(p-1)^2}$ contained in a class of subgroups isomorphic to $H$. Then, by the uniqueness of $G\leq H$, it follows that there exists a unique class of groups isomorphic to $H$ with signature $(0;2,2p,4)$ in $\mathcal{M}_{(p-1)^2}$.
\end{example}

Observe that the generating vector in Example \ref{ex-acc} is highly symmetric
and this symmetry is reflected by the fact that $G$ is contained in a larger
subgroup of $\mathcal{M}_{(p-1)^{2}}$. It seems
plausible that if a generating vector is highly symmetric, then the
corresponding conjugacy class of subgroups of the MCG will be contained in a
class of larger finite subgroups of the MCG. This could lead to a new way to
determine and enumerate classes of subgroups of the MCG which are not abelian
using the methods we have developed for abelian subgroups.

\end{document}